\documentclass[11pt]{article}
\usepackage{amsfonts}
\usepackage{latexsym,amsmath}
\usepackage{amssymb,array}
\usepackage{hyperref}
 \makeatletter \oddsidemargin 0in \evensidemargin 0 in \textwidth 16cm 
\begin{document}
\newcommand{\mf}{\mathbf}
\newcommand{\ov}{\overline}
\newcommand{\om}{\omega}
\newcommand{\ga}{\gamma}
\newcommand{\cd}{\circledast}
\newcommand{\eq}{\Leftrightarrow}
\newtheorem{t1}{Theorem}
\newtheorem{pro}{Proposition}
\newtheorem{remark}{Remark}
\newtheorem{counterexample}{Counterexample}
\newtheorem{coro}{Corollary}
\newtheorem{d1}{Definition}
\newtheorem{n1}{Notation}
\newtheorem{example}{Example}
\newtheorem{l1}{Lemma}
%%---------------------------------  Title Page  ----------------------------------------
\title{\bf Stochastic Comparisons of Lifetimes of Two Series and Parallel Systems with Location-Scale Family Distributed Components having Archimedean Copulas} 
 \author{Amarjit Kundu\\Department of
Mathematics\\
Santipur College\\ West Bengal, India\and Shovan Chowdhury\footnote{Corresponding
author e-mail: shovanc@iimk.ac.in; meetshovan@gmail.com}\;\\Quantitative Methods and Operations Management Area\\Indian Institute of Management, Kozhikode\\Kerala, India.
}
 %\date{Submitted: April, 2014\\$1$st Revision: December, 2014\\$2$nd Revision: August, 2015}
\maketitle
%----------------------------------------------Abstract-----------------------------------
{\bf Abstract-} In this paper, we compare the lifetimes of two series and two parallel systems stochastically where the lifetime of each component follows location-scale ($LS$) family of distributions. The comparison is carried out under two scenarios: one, that the components of the systems have a dependent structure sharing Archimedean copula and two, that the components are independently distributed. It is shown that the systems with components in series or parallel sharing Archimedean copula with more dispersion in the location or scale parameters results in better performance in the sense of the usual stochastic order. It is also shown that if the components are independently distributed, it is possible to obtain more generalized results as compared to the dependent set-up. The results in this paper generalizes similar results in both independent and dependent set up for exponential and Weibull distributed components.
%--------------------------------------------------Keywords-------------------------------
\newpage
%--------------  Notation-----------------------------------------------------------------

%----------------------------------------------------Introduction-------------------------
\section{Introduction}
\setcounter{equation}{0}
\hspace*{0.3 in} Stochastic comparison of system lifetimes has always been a relevant topic in reliability optimization and life testing experiments. These comparisons can be used to choose the best system structure under different criteria or to study where to place the different components in a system structure. If $X_{1:n}\leq X_{2:n}\leq\ldots\leq X_{n:n}$ denote the order statistics corresponding to the random variables $X_1, X_2,\ldots,X_n$, then the lifetime of a series and parallel system correspond to the smallest ($X_{1:n}$), and the largest ($X_{n:n}$) order statistic respectively. Classical theory of systems assumes that the lifetimes of the components are $iid$ (see David and Nagaraja~\cite{dn11}). Considerable amount of work has also been carried out in the past years in comparing the lifetimes of heterogeneous independent components of systems (largely on the smallest and the largest order statistics) with certain underlying distributions on both finite and infinite range with respect to usual stochastic ordering, hazard rate ordering, reversed hazard rate ordering and likelihood ratio ordering. One may refer to Dykstra \emph{et al.}~\cite{dkr11}, Zhao and Balakrishnan~(\cite{zb11.2}), Balakrishnan \emph{et al.}~\cite{ba1}, Torrado and Kochar~\cite{tr11}, Torrado~\cite{to}, Fang and Balakrishnan~\cite{fb}, Kundu \emph{et al.}~\cite{kun1}, Kundu and Chowdhury~(\cite{kun2},\cite{kun3}), Chowdhury and Kundu~\cite{ch} and Hazra \emph{et al.}~\cite{ha} for more detail. There are few work in the same area where the authors have compared systems stochastically through relative ageing, also known as ageing faster ordering in terms of hazard rate or reversed hazard rate ordering. One may refer to Sengupta and Deshpande~\cite{se}, Rezaei \emph{et al.}~\cite{re1} and Li and Li~\cite{lili} for more detail. 
  \\\hspace*{0.3 in} However, in practical situations, the components of a system may have a structural dependence which result in a set of statistically dependent observations. The dependence structure of the components are investigated by researchers very recently with the help of copulas. Navarro and Spizzichino~\cite{no} studied stochastic orders of series and parallel systems with components sharing a common copula. Rezapour and Alamatsaz~\cite{re} investigated stochastic orders on order statistics from samples with different survival Archimedean copulas. Li and Li~\cite{li} studied stochastic ordering of the sample minimums of Weibull samples sharing a common Archimedean survival copula. Li and Fang~\cite{li1} compared the lifetimes of parallel systems with proportional hazard rate (PHR) components following the Archimedean copula which was further investigated by Li \emph{et al.}~\cite{li2} and Fang \emph{et al.}~\cite{li3}. 
	\\\hspace*{0.3 in} The location-scale family of distributions is commonly used in lifetime studies. The most widely used statistical distributions are either members of this class or closely related to this class of distributions; such as exponential, normal, Weibull, lognormal, loglogistic, logistic, and extreme value distributions. Methods of inference and statistical theory for the general family can be applied to this large, important class of models. A random variable $X$ is said to follow location-scale family distribution, written as $LS$($\lambda,\sigma,F$) and will be termed as $LS$ family hereafter, if the distribution function of $X$ is given by
\begin{equation}\label{e0}
F(x;\lambda,\sigma,F)=F\left(\frac{x-\lambda}{\sigma}\right),\;\;x>\lambda,\;\sigma>0\;x,\;\lambda,\;\sigma\in\mathbb{R},
\end{equation}
where $\lambda$ and $\sigma$ are the location and the scale parameter respectively and $F\left(\cdot\right)$ is the
baseline distribution function of the rv $X$. Although significant previous research has compared series or parallel systems of heterogenous independent components including scale family of distributions (see Khaledi \emph{et al.}~\cite{kh}, Li \emph{et al.}~\cite{li2}, Kochar and Torrado~\cite{ko} and Li and Li~\cite{lili}), there has been few work examining similar comparisons for dependent components; furthermore, all such comparisons for dependent components assume either PHR components or Weibull distributed components. As $LS$ family of distributions cover a large pool of lifetime distributions, we are motivated to assume the component lifetimes to follow $LS$ family and to compare the lifetimes of two series or parallel systems stochastically, assuming a dependence structure in the components. In this sense, the paper distinguishes itself from the other few existing work. It generalizes the results on stochastic comparison of lifetimes of two series or parallel systems with heterogeneous dependent and independent components. Moreover, the comparisons are carried out under different baseline distributions of $LS$ family with stochastic ordering, hazard rate orderings and reversed hazard rate orderings between them. The rest of the paper is organized as follows. In Section 2, we have given the required definitions and some useful lemmas which are used throughout the paper. Results related to stochastic comparison of series systems with heterogeneous independent and dependent components are derived in Section~3. Section~4 discusses some results on parallel systems with dependent components.\\
%----------------------------------------------------Preliminaries---------------------------------------------------
\section{Preliminaries}
\setcounter{equation}{0}
\hspace*{0.3 in} For two absolutely continuous random variables $X$ and $Y$ with distribution functions $F\left(\cdot\right)$ and $G\left(\cdot\right)$, survival functions $\overline F\left(\cdot\right)$ and $\overline G\left(\cdot\right)$, density functions $f\left(\cdot\right)$ and $g\left(\cdot\right)$, hazard rate functions $r\left(\cdot\right)$ and $s\left(\cdot\right)$ and reversed hazard rate functions ${\tilde r(\cdot)}$ and ${\tilde s(\cdot)}$ respectively, $X$ is said to be smaller than $Y$ in $i)$ {\it likelihood ratio order} (denoted as $X\leq_{lr}Y$), if, for all $t$, $\frac{g(t)}{f(t)}$ increases in $t$, $ii)$ {\it hazard rate order} (denoted as $X\leq_{hr}Y$), if, for all $t$, $\frac{\overline G(t)}{\overline F(t)}$ increases in $t$ or equivalently $r(t)\ge s(t)$, $iii)$ {\it reversed hazard rate order} (denoted as $X\leq_{rh}Y$), if, for all $t$, $\frac{G(t)}{ F(t)}$ increases in $t$ or equivalently ${\tilde r(t)}\leq {\tilde s(t)}$, $iii)$ {\it Ageing faster order in terms of the hazard rate order} (denoted as $X\leq_{R-hr}Y$), if, for all $t$, $\frac{r_X(t)}{r_Y(t)}$ increases in $t$, and $iv)$  {\it usual stochastic order} (denoted as $X\leq_{st}Y$), if $F(t)\ge G(t)$ for all $t$. In the following diagram we present a chain of implications of the stochastic orders. For more on stochastic orders, see Shaked and Shanthikumar \cite{shak1}. 
\vspace{0.17 in}
\\\hspace*{1.7 in}$~~~~~~X\leq_{hr}Y$
\\\hspace*{1.7 in}$~~~~~~~~~~~\uparrow ~~~~~~~\searrow$
\\\hspace{6 in} $~~~~~~~~~~~~~~~~~~~~~~~~~~~~~~~~~~~~~~~~~X\leq_{lr}Y~~\rightarrow~~X\leq_{st}Y.$

\hspace{2.51 cm}$~~~~~~~~~~~~~~~~~~~~~~~~~\downarrow~~~~~~~~~\nearrow$

\hspace{2 cm}$~~~~~~~~~~~~~~~~~~~~~~~~X\leq_{rh}Y$
\\\hspace*{0.3 in} The notion of majorization (Marshall \emph{et al.}~\cite{Maol}) is essential for the understanding of the
stochastic inequalities for comparing order statistics. Let $\mathbb{R}^n$ be an $n$-dimensional Euclidean space. Further, for any two real vectors $\mathbf{x}=(x_1,x_2,\dots,x_n)\in \mathbb{R}^n$ and $\mathbf{y}=(y_1,y_2,\dots,y_n)\in \mathbb{R}^n$, write $x_{(1)}\le x_{(2)}\le\cdots\le x_{(n)}$ and $y_{(1)}\le y_{(2)}\le\cdots\le y_{(n)}$ as the increasing arrangements of the components of the vectors $\mathbf{x}$ and $\mathbf{y}$ respectively. The following definitions may be found in Marshall \emph{et al.}~\cite{Maol}.
%%%%%%%%%%%%%%%%%%%%%%%%%%%%%%%%%%%%%%%%%%%%%%%%%%%%%%%%%%%%%%%%%%%%%%%%%%%%%%%%%%%%%%%%%%%%%%%%%%%%%%%%%%%%%%%%%%%
\begin{d1}
\begin{enumerate}
\item[i)] The vector $\mathbf{x} $ is said to majorize the vector $\mathbf{y} $ (written as $\mathbf{x}\stackrel{m}{\succeq}\mathbf{y}$) if
$$\sum_{i=1}^j x_{(i)}\le\sum_{i=1}^j y_{(i)},\;j=1,\;2,\;\ldots, n-1,\;\;and \;\;\sum_{i=1}^n x_{(i)}=\sum_{i=1}^n y_{(i)}.$$
\item [ii)] The vector $\mathbf{x}$ is said to weakly supermajorize the vector $\mathbf{y}$
 (written as $\mathbf{x}\stackrel{\rm w}{\succeq} \mathbf{y}$) if
  $$\sum\limits_{i=1}^j x_{(i)}\leq \sum\limits_{i=1}^j y_{(i)}\quad \text{for}\;j=1,2,\dots,n.$$
 \item [iii)] The vector $\mathbf{x}$ is said to weakly submajorize the vector $\mathbf{y}$
 (written as $\mathbf{x}\;{\succeq}_{\rm w} \;\mathbf{y}$) if
  $$\sum\limits_{i=j}^n x_{(i)}\ge \sum\limits_{i=j}^n y_{(i)}\quad \text{for}\;j=1,2,\dots,n.$$
	\item [(iv)] The vector $\mathbf{x}$ is said to be $p$-larger than the vector $\mathbf{y}$
 (written as $\mathbf{x}\stackrel{\rm p}{\succeq} \mathbf{y}$) if
 $$\prod\limits_{i=1}^j x_{(i)}\leq \prod\limits_{i=1}^j y_{(i)}\quad \text{for}\;j=1,2,\dots,n.$$
 \item [(v)] The vector $\mathbf{x}$ is said to reciprocally majorize the vector $\mathbf{y}$
 (written as $\mathbf{x}\stackrel{\rm rm}\succeq \mathbf{y}$) if
 $$\sum\limits_{i=1}^j \frac{1}{x_{(i)}}\ge \sum\limits_{i=1}^j \frac{1}{y_{(i)}}\quad \text{for}\;j=1,2,\dots,n.$$
  \end{enumerate}
\end{d1}
\hspace*{0.3 in}It is not so difficult to show that 
$\mathbf{x}\stackrel{m}{\succeq}\mathbf{y}\Rightarrow\mathbf{x}\stackrel{ w}{\succeq} \mathbf{y}\Rightarrow\mathbf{x}\stackrel{ p}{\succeq} \mathbf{y}\Rightarrow\mathbf{x}\stackrel{ rm}{\succeq} \mathbf{y}$.
\begin{d1}\label{de2}
A function $\psi:\mathbb{R}^n\rightarrow\mathbb{R}$ is said to be Schur-convex (resp. Schur-concave) on $\mathbb{R}^n$ if 
\begin{equation*}
\mathbf{x}\stackrel{m}{\succeq}\mathbf{y} \;\text{implies}\;\psi\left(\mathbf{x}\right)\ge (\text{resp. }\le)\;\psi\left(\mathbf{y}\right)\;for\;all\;\mathbf{x},\;\mathbf{y}\in \mathbb{R}^n.
\end{equation*}
\end{d1}
\begin{d1}\label{de3}
For any integer $r$, a function $\psi:\mathbb{R}\rightarrow\mathbb{R}$ is said to be r-convex (resp. r-concave) on $\mathbb{R}$ if $\frac{d^{r}\psi(x)}{dx^r}\geq (\leq ) 0~for\;all\;x\in \mathbb{R}$.
\end{d1}
\begin{n1}
Let us introduce the following notations.
\begin{enumerate}
\item[(i)] $\mathcal{D}_{+}=\left\{\left(x_{1},x_2,\ldots,x_{n}\right):x_{1}\ge x_2\ge\ldots\ge x_{n}> 0\right\}$.
\item[(ii)] $\mathcal{E}_{+}=\left\{\left(x_{1},x_2,\ldots,x_{n}\right):0< x_{1}\leq x_2\leq\ldots\leq x_{n}\right\}$.
\end{enumerate}
\end{n1}
Let us first introduce the following lemmas which will be used in the next sections to prove the results. 
 \begin{l1}\label{l6}
 $\left(\text{Lemma 3.1 of Kundu \emph{et al.}~\cite{kun1}}\right)$~Let $\varphi:\mathcal{D_{+}}\rightarrow \mathbb{R}$ be a function, continuously differentiable on the interior of $\mathcal{D_{+}}$. Then, for $\mf{x},\mf{y}\in \mathcal{D_{+}}$,
 \begin{eqnarray*}
  \mf{x}\stackrel{m}\succeq\mf{y}\;\text{implies}\;\varphi(\mf{x})\ge (\text{resp.}\;\leq)\;\varphi(\mf{y})
 \end{eqnarray*}
if, and only if,
$$\varphi_{(k)}(\mf{z})\;\text{is decreasing (resp. increasing) in}\;k=1,2,\dots,n,$$
where $\varphi_{(k)}(\mf{z})=\partial\varphi(\mf{z})/\partial z_k$ denotes the partial derivative of $\varphi$ with respect to its $k$th argument.
\end{l1}
\begin{l1}\label{l7}
 $\left(\text{Lemma 3.3 of Kundu \emph{et al.}~\cite{kun1}}\right)\;$~Let $\varphi:\mathcal{E_{+}}\rightarrow \mathbb{R}$ be a function, continuously differentiable on the interior of $\mathcal{E_{+}}$.
 Then, for $\mf{x},\mf{y}\in \mathcal{E_{+}}$,
 \begin{eqnarray*}
  \mf{x}\stackrel{m}\succeq\mf{y}\;\text{implies}\;\varphi(\mf{x})\ge (\text{resp.}\;\leq)\;\varphi(\mf{y})
 \end{eqnarray*}
if, and only if,
$$\varphi_{(k)}(\mf{z})\;\text{is increasing (resp. decreasing) in}\;k=1,2,\dots,n,$$
where $\varphi_{(k)}(\mf{z})=\partial\varphi(\mf{z})/\partial z_k$ denotes the partial derivative of $\varphi$ with respect to its $k$th argument.
\end{l1}
\begin{l1}\label{l8}
 $\left(\text{Lemma 3.1 of Khaledi and Kochar}~\cite{kk}\right)\;$~Let $S\subseteq \mathbb{R}^n_+$. Further, let $\psi :S\rightarrow \mathbb{R}$ be a function. Then, for $\mbox{{\bf x}},\mbox{{\bf y}}\in S$,
 $${\bf x}\stackrel{ p}\succeq {\bf y}\;\text{implies}\;\psi({\bf x})\ge \;(resp. \leq)\;\psi({\bf y})$$
if, and only if,
\begin{enumerate}
 \item[(i)]$\psi(e^{a_1},\dots,e^{a_n})$ is Schur-convex (resp. Schur-concave) in $(a_1,\dots,a_n)\in S$,
\item[(ii)]$\psi(e^{a_1},\dots,e^{a_n})$ is decreasing (resp. increasing) in $a_i,$ for $i=1,\dots,n,$
\end{enumerate} 
where $a_i=\ln x_i$, for $i=1,\dots,n.$ 
\end{l1}
\begin{l1}\label{l9}
 $\left(\text{Lemma 4.1 of Hazra \emph{et al.}~\cite{ha}}\right)\;$~Let $S\subseteq \mathbb{R}^n_+$. Further, let $\psi :S\rightarrow \mathbb{R}$ be a function. Then, for $\mbox{{\bf x}},\mbox{{\bf y}}\in S$,
 $${\bf x}\stackrel{ rm}\succeq {\bf y}\;\text{implies}\;\psi({\bf x})\ge \;(resp. \leq) \psi({\bf y})$$
if, and only if,
\begin{enumerate}
 \item[(i)]$\psi(\frac{1}{a_1},\dots,\frac{1}{a_n})$ is Schur-convex (resp. Schur-concave) in $(a_1,\dots,a_n)\in S$,
\item[(ii)]$\psi(\frac{1}{a_1},\dots,\frac{1}{a_n})$ is increasing (resp. decreasing) in $a_i,$ for $i=1,\dots,n,$
\end{enumerate} 
where $a_i=\frac{1}{x_i}$, for $i=1,\dots,n.$  
\end{l1}
\begin{l1}\label{l10}
$\left(\text{Theorem A.8 of Marshall \emph{et al.}~\cite{Maol}}\;p.p.~87\right)$~Let $S\subseteq \mathbb{R}^n$. Further, let $\varphi:S\rightarrow\mathbb{R}$ be a function. Then for $\mbox{{\bf x}}$, $\mbox{{\bf y}}\in S$,
$$\mbox{{\bf x}}\succeq_{w} \mbox{{\bf y}}\Longrightarrow\varphi\left(\mbox{{\bf x}}\right)\ge(resp. \le)\varphi\left(\mbox{{\bf y}}\right)$$
if, and if, $\varphi$ is both increasing (resp. decreasing) and Schur-convex (resp. Schur-concave) on $S$. Similarly,
$${\bf x}\stackrel{w}\succeq {\bf y}\Longrightarrow\varphi\left({\bf x}\right)\ge(resp. \le)\varphi\left({\bf y}\right)$$
if, and if, $\varphi$ is both decreasing (resp. increasing) and Schur-convex (resp. Schur-concave) on $S$.
\end{l1}
\hspace*{0.3 in} Now, let us recall that a copula associated with a multivariate distribution function $F$ is a function $C:\left[0,1\right]^n\longmapsto\left[0,1\right]$ satisfying: $F(x)=C\left(F_{1}(X_1),..., F_{n}(X_n)\right),$ where the $F_i$'s, $1\leq i\leq n$ are the univariate marginal distribution functions of $X_i$s. Similarly, a survival copula associated with a multivariate survival function $\overline{F}$ is a function $\overline{C}:\left[0,1\right]^n\longmapsto\left[0,1\right]$ satisfying:
$$\overline{F}(x)=P\left(X_1>x_1,...,X_n>x_n\right)=\overline{C}\left(\overline{F}_1(x_1),...,\overline{F}_n(x_n)\right),$$ 
where, for $1\leq i\leq n$, $\overline{F}_i(\cdot)=1-F_i(\cdot)$  are the univariate survival functions. In particular, a copula $C$ is Archimedean if there exists a generator $\psi:\left[0,\infty\right]\longmapsto\left[0,1\right]$ such that
$$C\left(\mathbf{u}\right)=\psi\left(\psi^{-1}(u_1),...,\psi^{-1}(u_d)\right).$$
For $C$ to be Archimedean copula, it is sufficient and necessary that $\psi$ satisfies $i)$ $\psi(0)=1$ and $\psi(\infty)=0$ and $ii)$ $\psi$ is $d-$monotone, i.e. $\frac{(-1)^k d^k \psi(s)}{ds^k}\ge 0$ for $k\in \left\{0,1,...,d-2\right\}$ and $\frac{(-1)^{d-2} d^{d-2} \psi(s)}{ds^{d-2}}$ is decreasing and convex. Archimedean copulas
cover a wide range of dependence structures including the independence copula and the Clayton copula. For more detail on Archimedean copula, see, Nelsen~\cite{ne} and McNeil and N$\check{e}$slehov$\acute{a}$~\cite{mc}. In this paper, Archimedean copula is specifically employed to model on the dependence structure among random variables in a sample. The following important lemma is used in the next sections to prove some of the important theorems. 
\begin{l1}\label{l11}
 $\left(\text{Li and Fang~\cite{li1}}\right)\;$~For two n-dimensional Archimedean copulas $C_{\psi_1}\left(\mathbf{u}\right)$ and $C_{\psi_2}\left(\mathbf{u}\right)$, with $\phi_2=\psi_{2}^{-1}=sup\left\{x\in\mathbb{R}:\psi(x)>u\right\}$, the right continuous inverse, if $\phi_2\circ\psi_1$ is super-additive, then $C_{\psi_1}\left(\mathbf{u}\right)\leq C_{\psi_2}\left(\mathbf{u}\right)$ for all $\mathbf{u}\in [0,1]^n.$ Recall that a function $f$ is said to be super-additive if $f(x+y)\ge f(x) + f(y)$, for all $x$ and $y$ in the domain of $f$.
 \end{l1} 
%%--------------------------comparison of series system (independent and dependet)----------------------------------------------------------

\section{Comparison of Series Systems with LS Distributed Components}
\setcounter{equation}{0}
\hspace*{0.3 in} This section is devoted to the comparison of two series systems with heterogenous $LS$ family distributed components. The comparison is carried out under two scenarios: one, that the components have a dependent structure sharing Archimedean copulas and the other is that the components are  independently distributed. \\
 %%--------------------------------Subsection 3.1-------------------------------------------------------
\subsection{\small Some Results on Heterogenous Dependent Components}
\hspace*{0.3 in} Let, $X$ and $Y$ be two random variables having distribution functions $F(\cdot)$ and $G(\cdot)$ respectively. Also suppose that $X_i\sim LS\left(\lambda_i,\sigma_i,F\right)$ and $Y_i\sim LS\left(\mu_i,\xi_i,G\right)$ ($i=1,2,\ldots,n$) be two sets of $n$ dependent random variables with Archimedean copulas having generators $\psi_1~\left(\phi_1=\psi_{1}^{-1}\right)$ and $\psi_2~\left(\phi_2=\psi_{2}^{-1}\right)$ respectively. Also suppose that $\overline{G}_{1:n}\left(\cdot\right)$ and $\overline{H}_{1:n}\left(\cdot\right)$ be the survival functions of $X_{1:n}$ and $Y_{1:n}$ respectively. Then, 
\begin{equation*}
\overline{G}_{1:n}\left(t\right)=\psi_1\left[\sum_{k=1}^n \phi_1\left\{1-F\left(\frac{t-\lambda_k}{\sigma_k}\right)\right\}\right],~t>max(\lambda_k,~\forall k),
\end{equation*}
and
\begin{equation*}
\overline{H}_{1:n}\left(t\right)=\psi_2\left[\sum_{k=1}^n \phi_2\left\{1-G\left(\frac{t-\mu_k}{\xi_k}\right)\right\}\right],~t>max(\mu_k,~\forall k).
\end{equation*}
\hspace*{0.3 in} Let $r_X(u)$ and $r_Y(u)$ are the hazard rate functions of the random variables $X$ and $Y$ respectively. The first two theorems show that usual stochastic ordering exists between $X_{1:n}$ and $Y_{1:n}$ under weak majorization order of the scale parameters and stochastic ordering between $X$ and $Y$. 
\begin{t1}\label{th1}
Let $X_1,X_2,...,X_n$ be a set of dependent random variables sharing Archimedean copula having generator $\psi_1$ such that $X_i\sim $LS$\left(\lambda_i,\sigma_i,F\right),~i=1,2,...,n$. Let $Y_1,Y_2,...,Y_n$ be another set of dependent random variables sharing Archimedean copula having generator $\psi_2$ such that $Y_i\sim $LS$\left(\lambda_i,\xi_i,G\right),~i=1,2,...,n$. Assume that $\mbox{\boldmath $\sigma$},\mbox{\boldmath $\xi$}\;and\; \mbox{\boldmath $\lambda$}\in \mathcal{D}_+$ (or $\mathcal{E}_+$). Further suppose that $\phi_2\circ\psi_1$ is super-additive, $\psi_1$ or $\psi_2$ is log-concave and $X\leq_{st}Y$. If either $r_X(u)$ or $r_Y(u)$ is decreasing in $u$, then $\mbox{\boldmath $\frac{1}{\sigma}$}\stackrel{w}{\preceq} \mbox{\boldmath $\frac{1}{\xi}$}$ implies $X_{1:n}\leq_{st}Y_{1:n}$.
\end{t1}
{\bf Proof:} By Lemma \ref{l11}, super-additivity of $\phi_2\circ\psi_1$ implies that 
\begin{equation}\label{e1}
\psi_1\left[\sum_{k=1}^n \phi_1\left\{1-F\left(\frac{t-\lambda_k}{\sigma_k}\right)\right\}\right]\leq \psi_2\left[\sum_{k=1}^n \phi_2\left\{1-F\left(\frac{t-\lambda_k}{\sigma_k}\right)\right\}\right].
\end{equation}
As $X\leq_{st}Y$ and $\phi_2$ and $\psi_2$ are decreasing in $x$, it can be easily shown that 
\begin{equation}\label{e2}
\psi_2\left[\sum_{k=1}^n \phi_2\left\{1-F\left(\frac{t-\lambda_k}{\sigma_k}\right)\right\}\right]\leq \psi_2\left[\sum_{k=1}^n \phi_2\left\{1-G\left(\frac{t-\lambda_k}{\sigma_k}\right)\right\}\right].
\end{equation}
(\ref{e1}) and (\ref{e2}) together implies that 
\begin{equation*}\label{e3}
\psi_1\left[\sum_{k=1}^n \phi_1\left\{1-F\left(\frac{t-\lambda_k}{\sigma_k}\right)\right\}\right]\leq \psi_2\left[\sum_{k=1}^n \phi_2\left\{1-G\left(\frac{t-\lambda_k}{\sigma_k}\right)\right\}\right].
\end{equation*}
Therefore, to prove the result it suffices to prove that $$\psi_2\left[\sum_{k=1}^n \phi_2\left\{1-G\left(\frac{t-\lambda_k}{\sigma_k}\right)\right\}\right]\leq \psi_2\left[\sum_{k=1}^n \phi_2\left\{1-G\left(\frac{t-\lambda_k}{\xi_k}\right)\right\}\right].$$ 
Let us assume that
$$\psi_2\left[\sum_{k=1}^n \phi_2\left\{1-G\left(\frac{t-\lambda_k}{\sigma_k}\right)\right\}\right]=\psi_2\left[\sum_{k=1}^n \phi_2\left\{1-G\left(p_{k}(t-\lambda_k)\right)\right\}\right]=\Psi(\mbox{\boldmath $p$}),$$ 
where $\mbox{\boldmath $p$}=\left(p_1,p_2,...,p_n\right)=\left(\frac{1}{\sigma_1},\frac{1}{\sigma_2},...,\frac{1}{\sigma_n}\right).$ Hence, by Lemma \ref{l10}, to prove the result, it suffices to prove that $\Psi(\mbox{\boldmath $p$})$ is decreasing and s-convex in $\mbox{\boldmath $p$}$. Now, if $\mbox{\boldmath $\sigma$}, \mbox{\boldmath $\lambda$}\in \mathcal{D}_+\left(or~\mathcal{E}_+\right)$ and $r_Y(u)$ is decreasing in $u$, then, for all $i\leq j$,
$$r_Y\left(p_{i}(t-\lambda_i)\right)\ge (\leq) r_Y\left(p_{j}(t-\lambda_j)\right).$$ 
Again, as $\psi_2$ is log-concave, giving that $\frac{\psi_{2}(u)}{\psi_{2}^{'}(u)}$ is increasing in $u,$ then, for all $i\leq j$, it can be written that
\begin{equation}\label{e4}
 -r_Y\left(p_{i}(t-\lambda_i)\right)\frac{\psi_2\left(u_i\right)}{\psi_{2}^{'}\left(u_i\right)}\ge (\leq) -r_Y\left(p_{j}(t-\lambda_j)\right)\frac{\psi_2\left(u_j\right)}{\psi_{2}^{'}\left(u_j\right)},  
\end{equation} 
where $u_i=\phi_2\left[1-G\left(p_{i}(t-\lambda_i)\right)\right].$ Now, differentiating $\Psi(\mbox{\boldmath $p$})$ with respect to $p_i$, we get
\begin{equation*}
\begin{split}
\frac{\partial \Psi}{\partial p_i}&=-\psi_{2}^{'}\left[\sum_{k=1}^n \phi_2\left\{1-G\left(p_k(t-\lambda_k)\right)\right\}\right](t-\lambda_i) r_Y\left(p_{i}(t-\lambda_i)\right)\\&\quad\frac{\psi_{2}\left[\phi_2\left\{1-G\left(p_{i}(t-\lambda_i)\right)\right\}\right]}{\psi_{2}^{'}\left[\phi_2\left\{1-G\left(p_{i}(t-\lambda_i)\right)\right\}\right]}\leq 0,
\end{split}
\end{equation*}
proving that $\Psi(\mbox{\boldmath $p$})$ is decreasing in each $p_i.$ Moreover, using (\ref{e4}), it can be easily shown that
$\frac{\partial \Psi}{\partial p_i}-\frac{\partial \Psi}{\partial p_j}\leq (\ge) 0$. Thus, by Lemma \ref{l7} (Lemma \ref{l6}) it can be written that $\Psi(\mbox{\boldmath $p$})$ is s-convex in $\mbox{\boldmath $p$}$. This proves the result. \hfill$\Box$\\
The next theorem discusses about the stochastic ordering between $X_{1:n}$ and $Y_{1:n}$ under weak majorization order of the scale parameters when $\psi_1$ or $\psi_2$ is log-convex. The theorem can be proved in the similar line as of the previous one and hence the proof is omitted. 
\begin{t1}\label{th1a}
Let $X_1,X_2,...,X_n$ be a set of dependent random variables sharing Archimedean copula having generator $\psi_1$ such that $X_i\sim $LS$\left(\lambda_i,\sigma_i,F\right),~i=1,2,...,n$. Let $Y_1,Y_2,...,Y_n$ be another set of dependent random variables sharing Archimedean copula having generator $\psi_2$ such that $Y_i\sim $LS$\left(\lambda_i,\xi_i,G\right),~i=1,2,...,n$. Assume that $\mbox{\boldmath $\sigma$},\mbox{\boldmath $\xi$}$ and $\mbox{\boldmath $\lambda$}\in \mathcal{D}_+$ (or $\mathcal{E}_+$). Further suppose that $\phi_2\circ\psi_1$ is super-additive, $\psi_1$ or $\psi_2$ is log-convex and $X\leq_{st}Y$. If either $r_X(u)$ or $r_Y(u)$ is increasing in $u$, then $\mbox{\boldmath $\frac{1}{\sigma}$}{\succeq}_w \mbox{\boldmath $\frac{1}{\xi}$}$ implies $X_{1:n}\leq_{st}Y_{1:n}$.
\end{t1}
\hspace*{0.3 in} In the next theorem usual stochastic ordering between $X_{1:n}$ and $Y_{1:n}$ has been established under $p$ and $rm$ orderings of the scale parameters.
\begin{t1}\label{th2}
Let $X_1,X_2,...,X_n$ be a set of dependent random variables sharing Archimedean copula having generator $\psi_1$ such that $X_i\sim $LS$\left(\lambda_i,\sigma_i,F\right),~i=1,2,...,n$. Let $Y_1,Y_2,...,Y_n$ be another set of dependent random variables sharing Archimedean copula having generator $\psi_2$ such that $Y_i\sim $LS$\left(\lambda_i,\xi_i,G\right),~i=1,2,...,n$. Assume that $\mbox{\boldmath $\sigma$},\mbox{\boldmath $\xi$},\mbox{\boldmath $\lambda$}\in \mathcal{D}_+$ (or $\mathcal{E}_+$). Further suppose that $\phi_2\circ\psi_1$ is super-additive, $\psi_1$ or $\psi_2$ is log-concave and $X\leq_{st}Y$, then,
\begin{enumerate}
\item[i)] $\mbox{\boldmath $\frac{1}{\sigma}$}\stackrel{p}{\preceq} \mbox{\boldmath $\frac{1}{\xi}$}$ $\Rightarrow$ $X_{1:n}\leq_{st}Y_{1:n}$,$~$ if $ur_X(u)$ or $ur_Y(u)$ is decreasing in $u$;
\item[ii)] $\mbox{\boldmath $\frac{1}{\sigma}$}\stackrel{rm}{\preceq} \mbox{\boldmath $\frac{1}{\xi}$}$ $\Rightarrow$ $X_{1:n}\leq_{st}Y_{1:n}$,$~$ if $u^2r_X(u)$ or $u^2r_Y(u)$ is decreasing in $u$.
\end{enumerate}
\end{t1}
{\bf Proof:} i) As $\phi_2\circ\psi_1$ is super-additive, by Lemma \ref{l11} it can be written that 
\begin{equation}\label{e5}
\psi_1\left[\sum_{k=1}^n \phi_1\left\{1-F\left(\frac{t-\lambda_k}{\sigma_k}\right)\right\}\right]\leq \psi_2\left[\sum_{k=1}^n \phi_2\left\{1-F\left(\frac{t-\lambda_k}{\sigma_k}\right)\right\}\right].
\end{equation}
Again noticing the fact that $X\leq_{st}Y$ and $\phi_2$ and $\psi_2$ are decreasing function of $x$, it can be easily shown that 
\begin{equation}\label{e6}
\psi_2\left[\sum_{k=1}^n \phi_2\left\{1-F\left(\frac{t-\lambda_k}{\sigma_k}\right)\right\}\right]\leq \psi_2\left[\sum_{k=1}^n \phi_2\left\{1-G\left(\frac{t-\lambda_k}{\sigma_k}\right)\right\}\right].
\end{equation}
So, (\ref{e5}) and (\ref{e6}) together gives,
\begin{equation*}\label{e7}
\psi_1\left[\sum_{k=1}^n \phi_1\left\{1-F\left(\frac{t-\lambda_k}{\sigma_k}\right)\right\}\right]\leq \psi_2\left[\sum_{k=1}^n \phi_2\left\{1-G\left(\frac{t-\lambda_k}{\sigma_k}\right)\right\}\right].
\end{equation*}
Therefore, to prove the result, it suffices to prove that $$\psi_2\left[\sum_{k=1}^n \phi_2\left\{1-G\left(\frac{t-\lambda_k}{\sigma_k}\right)\right\}\right]\leq \psi_2\left[\sum_{k=1}^n \phi_2\left\{1-G\left(\frac{t-\lambda_k}{\xi_k}\right)\right\}\right].$$ 
Let us assume that\\
$\psi_2\left[\sum_{k=1}^n \phi_2\left\{1-G\left(\frac{t-\lambda_k}{\sigma_k}\right)\right\}\right]=\psi_2\left[\sum_{k=1}^n \phi_2\left\{1-G\left(e^{p_k}(t-\lambda_k)\right)\right\}\right]=\Psi(\mbox{\boldmath $p$})$, where $\mbox{\boldmath $p$}=\left(p_1,p_2,...,p_n\right)=\left(\log\frac{1}{\sigma_1},\log\frac{1}{\sigma_2},...,\log\frac{1}{\sigma_n}\right).$ Hence, by Lemma \ref{l8}, it suffices to prove that $\Psi(\mbox{\boldmath $p$})$ is decreasing and s-convex in $\mbox{\boldmath $p$}$. Now, if $\mbox{\boldmath $\lambda$},\mbox{\boldmath $\sigma$}\in \mathcal{D}_+\left(or~\mathcal{E}_+\right)$ and $ur_Y(u)$ is decreasing in $u$, then for all $i\leq j$,
$$e^{p_i}(t-\lambda_i)r_Y\left(e^{p_i}(t-\lambda_i)\right)\ge (\leq)  e^{p_j}(t-\lambda_j)r_Y\left(e^{p_j}(t-\lambda_j)\right).$$ 
Again, if $\psi_2$ is log-concave, then $\frac{\psi_{2}(u)}{\psi_{2}^{'}(u)}$ is increasing in $u,$ which implies that
\begin{equation}\label{e8}
 -e^{p_i}(t-\lambda_i)r_Y\left(e^{p_i}(t-\lambda_i)\right)\frac{\psi_2\left(u_i\right)}{\psi_{2}^{'}\left(u_i\right)}\ge (\leq) -e^{p_j}(t-\lambda_j)r_Y\left(e^{p_j}(t-\lambda_j)\right)\frac{\psi_2\left(u_j\right)}{\psi_{2}^{'}\left(u_j\right)};~i\leq j,  
\end{equation} 
with $u_i=\phi_2\left[1-G\left(e^{p_i}(t-\lambda_i)\right)\right].$ Now, differentiating $\Psi(\mbox{\boldmath $p$})$ with respect to $p_i$, we get,  
\begin{equation*}
\begin{split}
\frac{\partial \Psi}{\partial p_i}&=-\psi_{2}^{'}\left[\sum_{k=1}^n \phi_2\left\{1-G\left(e^{p_k}(t-\lambda_k)\right)\right\}\right]e^{p_i}(t-\lambda_i)r_Y\left(e^{p_i}(t-\lambda_i)\right)\\&\quad\frac{\psi_{2}\left[\phi_2\left\{1-G\left(e^{p_i}(t-\lambda_i)\right)\right\}\right]}{\psi_{2}^{'}\left[\phi_2\left\{1-G\left(e^{p_i}(t-\lambda_i)\right)\right\}\right]}\leq 0,
\end{split}
\end{equation*} 
proving that $\Psi(\mbox{\boldmath $p$})$ is decreasing in each $p_i.$ Moreover, using (\ref{e8}), it can be easily shown that
$$\frac{\partial \Psi}{\partial p_i}-\frac{\partial \Psi}{\partial p_j}\leq (\ge) 0.$$ Thus, by Lemma \ref{l7} (Lemma \ref{l6}), $\Psi(\mbox{\boldmath $p$})$ is s-convex in $\mbox{\boldmath $p$}$. This proves the result. 
\\\hspace*{0.3 in} To prove ii), using Lemma \ref{l9},  it is to prove that 
$$\psi_2\left[\sum_{k=1}^n \phi_2\left\{1-G\left(\frac{t-\lambda_k}{\sigma_k}\right)\right\}\right]=\Psi(\mbox{\boldmath $\sigma$}),$$ 
is increasing in each $\sigma_i$ and s-convex in $\mbox{\boldmath $\sigma$}$. Now, if $\mbox{\boldmath $\lambda$},\mbox{\boldmath $\sigma$}\in \mathcal{D}_+\left( \mathcal{E}_+\right)$ and $u^{2}r_Y(u)$ is decreasing in $u$, following similar argument as in the previous result, it can be shown that, for all $i\leq j$,
\begin{equation}\label{e9}
 -\left(\frac{t-\lambda_i}{\sigma_i}\right)^{2} r_Y\left(\frac{t-\lambda_i}{\sigma_i}\right)\frac{1}{t-\lambda_i}\frac{\psi_2\left(u_i\right)}{\psi_{2}^{'}\left(u_i\right)}\ge (\leq) -\left(\frac{t-\lambda_j}{\sigma_j}\right)^{2} r_Y\left(\frac{t-\lambda_j}{\sigma_j}\right)\frac{1}{t-\lambda_j}\frac{\psi_2\left(u_j\right)}{\psi_{2}^{'}\left(u_j\right)},  
\end{equation} 
with $u_i=\phi_2\left[1-G\left(\frac{t-\lambda_i}{\sigma_i}\right)\right].$ Now, differentiating $\Psi(\mbox{\boldmath $\sigma$})$ with respect to $\sigma_i$, we get,
\begin{equation*}
\begin{split}
\frac{\partial \Psi}{\partial \sigma_i}&=\psi_{2}^{'}\left[\sum_{k=1}^n \phi_2\left\{1-G\left(\frac{t-\lambda_k}{\sigma_k}\right)\right\}\right]\left(\frac{t-\lambda_i}{\sigma_i}\right)^{2} r_Y\left(\frac{t-\lambda_i}{\sigma_i}\right)\frac{1}{t-\lambda_i}\\&\quad\frac{\psi_{2}\left[\phi_2\left(1-G\left(\frac{t-\lambda_i}{\sigma_i}\right)\right)\right]}{\psi_{2}^{'}\left[\phi_2\left(1-G\left(\frac{t-\lambda_i}{\sigma_i}\right)\right)\right]}\ge 0,
\end{split}
\end{equation*} 
proving that $\Psi(\mbox{\boldmath $\sigma$})$ is increasing in each $\sigma_i.$ Moreover, using (\ref{e9}), it can be easily shown that
$$\frac{\partial \Psi}{\partial \sigma_i}-\frac{\partial \Psi}{\partial \sigma_j}\ge (\leq) 0.$$ Thus, by Lemma \ref{l6} (Lemma \ref{l7}), $\Psi(\mbox{\boldmath $\sigma$})$ is s-convex in $\mbox{\boldmath $\sigma$}$. \hfill$\Box$\\
\hspace*{0.3 in} The next two theorems show that usual stochastic ordering exists between $X_{1:n}$ and $Y_{1:n}$ under majorization order of the location parameters. Proof of the second theorem follows from the first one, and hence is omitted.
\begin{t1}\label{th3}
Let $X_1,X_2,...,X_n$ be a set of dependent random variables sharing Archimedean copula having generator $\psi_1$ such that $X_i\sim $LS$\left(\lambda_i,\sigma_i,F\right),~i=1,2,...,n$. Let $Y_1,Y_2,...,Y_n$ be another set of dependent random variables sharing Archimedean copula having generator $\psi_2$ such that $Y_i\sim $LS$\left(\mu_i,\sigma_i,G\right),~i=1,2,...,n$. Assume that $\mbox{\boldmath $\lambda$},\mbox{\boldmath $\mu$},\mbox{\boldmath $\sigma$}\in \mathcal{D}_+$ (or $\mathcal{E}_+$). Further suppose that $\phi_2\circ\psi_1$ is super-additive, $\psi_1$ or $\psi_2$ is log-concave and either $ur_X(u)$ or $ur_Y(u)$ is decreasing in $u,$ then, $X\leq_{st}Y$ and $\mbox{\boldmath $\lambda$}{\preceq}_{w} \mbox{\boldmath $\mu$}$ $\Rightarrow$ $X_{1:n}\leq_{st}Y_{1:n}$.
\end{t1}
{\bf Proof:} Following the same argument as in (\ref{e5}) and (\ref{e6}) of Theorem~\ref{th2}, it can be shown that
\begin{equation}\label{e10}
\psi_1\left[\sum_{k=1}^n \phi_1\left\{1-F\left(\frac{t-\lambda_k}{\sigma_k}\right)\right\}\right]\leq \psi_2\left[\sum_{k=1}^n \phi_2\left\{1-G\left(\frac{t-\lambda_k}{\sigma_k}\right)\right\}\right].
\end{equation}
As before, to prove the result, it suffices to prove that $$\psi_2\left[\sum_{k=1}^n \phi_2\left\{1-G\left(\frac{t-\lambda_k}{\sigma_k}\right)\right\}\right]\leq \psi_2\left[\sum_{k=1}^n \phi_2\left\{1-G\left(\frac{t-\mu_k}{\sigma_k}\right)\right\}\right].$$ 
Using Lemma \ref{l10}, it can be said that the above relation will hold if 
$$\Psi(\mbox{\boldmath $\lambda$})=\psi_2\left[\sum_{k=1}^n \phi_2\left\{1-G\left(\frac{t-\lambda_k}{\sigma_k}\right)\right\}\right]$$
is increasing in each $\lambda_i$ and s-convex in $\mbox{\boldmath $\lambda$}$. Now, as $\mbox{\boldmath $\lambda$},\mbox{\boldmath $\sigma$}\in \mathcal{D}_+\left(or~\mathcal{E}_+\right)$, $ur_Y(u)$ is decreasing in $u$ and $\psi_2$ is log-concave, following the same argument as of Theorem \ref{th2}, we can show that 
\begin{equation}\label{e11}
 -\frac{1}{t-\lambda_i}\left(\frac{t-\lambda_i}{\sigma_i}\right) r_Y\left(\frac{t-\lambda_i}{\sigma_i}\right)\frac{\psi_2\left(u_i\right)}{\psi_{2}^{'}\left(u_i\right)}\ge (\leq)-\frac{1}{t-\lambda_j}\left(\frac{t-\lambda_j}{\sigma_j}\right) r_Y\left(\frac{t-\lambda_j}{\sigma_j}\right)\frac{\psi_2\left(u_j\right)}{\psi_{2}^{'}\left(u_j\right)};~i\leq j,  
\end{equation} 
with $u_i=\phi_2\left[1-G\left(\frac{t-\lambda_i}{\sigma_i}\right)\right].$ Differentiating $\Psi(\mbox{\boldmath $\lambda$})$ with respect to $\lambda_i$, we get
\begin{equation*}
\begin{split}
\frac{\partial \Psi}{\partial \lambda_i}&=\psi_{2}^{'}\left[\sum_{k=1}^n \phi_2\left\{1-G\left(\frac{t-\lambda_i}{\sigma_i}\right)\right\}\right]\frac{1}{t-\lambda_i}\frac{t-\lambda_i}{\sigma_i} r_Y\left(\frac{t-\lambda_i}{\sigma_i}\right)\\&\quad\frac{\psi_{2}\left[\phi_2\left\{1-G\left(\frac{t-\lambda_i}{\sigma_i}\right)\right\}\right]}{\psi_{2}^{'}\left[\phi_2\left\{1-G\left(\frac{t-\lambda_i}{\sigma_i}\right)\right\}\right]}\ge 0,
\end{split}
\end{equation*}
proving that $\Psi(\mbox{\boldmath $\lambda$})$ is increasing in each $\lambda_i.$ Moreover, using (\ref{e11}), it can be easily shown that
$$\frac{\partial \Psi}{\partial \lambda_i}-\frac{\partial \Psi}{\partial \lambda_j}\ge (\le) 0.$$ 
Thus, by Lemma \ref{l6} (Lemma \ref{l7}), $\Psi(\mbox{\boldmath $\lambda$})$ is s-convex in $\mbox{\boldmath $\lambda$}$. This proves the result. \hfill$\Box$
\begin{t1}\label{th3a}
Let $X_1,X_2,...,X_n$ be a set of dependent random variables sharing Archimedean copula having generator $\psi_1$ such that $X_i\sim $LS$\left(\lambda_i,\sigma_i,F\right),~i=1,2,...,n$. Let $Y_1,Y_2,...,Y_n$ be another set of dependent random variables sharing Archimedean copula having generator $\psi_2$ such that $Y_i\sim $LS$\left(\mu_i,\sigma_i,G\right),~i=1,2,...,n$. Assume that $\mbox{\boldmath $\lambda$},\mbox{\boldmath $\mu$},\mbox{\boldmath $\sigma$}\in \mathcal{D}_+$ (or $\mathcal{E}_+$). Further suppose that $\phi_2\circ\psi_1$ is super-additive, $\psi_1$ or $\psi_2$ is log-convex and either $r_X(u)$ or $r_Y(u)$ is increasing in $u,$ then, $X\leq_{st}Y$ and $\mbox{\boldmath $\lambda$}\stackrel{w}{\succeq}\mbox{\boldmath $\mu$}$ $\Rightarrow$ $X_{1:n}\leq_{st}Y_{1:n}$.
\end{t1}
\subsection{\small Heterogenous independent Components}
\hspace*{0.3 in} In this subsection, we will compare two series systems with heterogenous independent $LS$ family distributed components. For $i=1,2,\ldots,n$, let $X_i$ and $Y_i$ be two sets of $n$ independent random variables following $LS$ distribution with parameters ($\mbox{\boldmath$\lambda$}$, $\mbox{\boldmath$\sigma$}$) and ($\mbox{\boldmath$\mu$}$, $\mbox{\boldmath$\xi$}$) respectively, as given in (\ref{e0}).
If $\overline{F}_{1:n}\left(\cdot\right)$ and $\overline{G}_{1:n}\left(\cdot\right)$ be the survival functions of $X_{1:n}$ and $Y_{1:n}$ respectively, then clearly
\begin{equation*}
\overline{F}_{1:n}\left(x\right)=\prod_{k=1}^n \left[1-F\left(\frac{t-\lambda_k}{\sigma_k}\right)\right]
\end{equation*}
and
\begin{equation*}
\overline{G}_{1:n}\left(x\right)=\prod_{k=1}^n \left[1-G\left(\frac{t-\mu_k}{\xi_k}\right)\right].
\end{equation*}
Again, if $r_{1:n}(\cdot)$ and $s_{1:n}(\cdot)$ are the hazard rate functions of $X_{1:n}$ and $Y_{1:n}$ respectively then,
\begin{equation}\label{e20}
r_{1:n}\left(t\right)=\sum_{k=1}^n{\frac{1}{\sigma_k}\frac{f\left(\frac{t-\lambda_k}{\sigma_k}\right)}{1-F\left(\frac{t-\lambda_k}{\sigma_k}\right)}}=\sum_{k=1}^n{\frac{1}{\sigma_k}r_{X}\left(\frac{t-\lambda_k}{\sigma_k}\right)}
\end{equation} 
and
\begin{equation}\label{e21}
s_{1:n}\left(t\right)=\sum_{k=1}^n{\frac{1}{\xi_k}\frac{g\left(\frac{t-\mu_k}{\xi_k}\right)}{1-G\left(\frac{t-\mu_k}{\xi_k}\right)}}=\sum_{k=1}^n{\frac{1}{\xi_k}r_{Y}\left(\frac{t-\mu_k}{\xi_k}\right)},
\end{equation} 
%where $f(\cdot)$ and $g(\cdot)$ are the probability density functions of $X$ and $Y$ respectively. 
\hspace*{0.3 in} In the next two theorems, stochastic comparison between minimum order statistics from the $LS$ family with respect to hazard rate ordering has been discussed. These results strengthen Theorems~\ref{th1}~-~\ref{th3} for the independent case, as hazard rate ordering implies usual stochastic ordering.
\begin{t1}\label{th4}
Let $X_1,X_2,...,X_n$ be a set of independent random variables such that $X_i\sim $LS$\left(\lambda_i,\sigma_i,F\right)$, $i=1,2,...,n$. Let $Y_1,Y_2,...,Y_n$ be another set of independent random variables such that $Y_i\sim $LS$\left(\lambda_i,\xi_i,G\right)$, $i=1,2,...,n$. Suppose that $\mbox{\boldmath $\sigma$},\mbox{\boldmath $\xi$},\mbox{\boldmath $\lambda$}\in \mathcal{D}_+$ (or $\mathcal{E}_+$). Then, if $X\leq_{hr}Y$ and
\begin{enumerate}
\item[i)] either $ur_{X}(u)$ or $ur_{Y}(u)$ is concave in $u,$ then, $\mbox{\boldmath $\frac{1}{\sigma}$}\stackrel{m}{\preceq} \mbox{\boldmath $\frac{1}{\xi}$}$ $\Rightarrow$ $X_{1:n}\leq_{hr}Y_{1:n}$,
\item[ii)] either $ur_{X}(u)$ or $ur_{Y}(u)$ is increasing (decreasing) and concave in $u,$ then, $\mbox{\boldmath $\frac{1}{\sigma}$}\stackrel{w}{\preceq}(\preceq_w) \mbox{\boldmath $\frac{1}{\xi}$}$ $\Rightarrow$ $X_{1:n}\leq_{hr}Y_{1:n}$,
\item[iii)] either $ur_{X}(u)$ or $ur_{Y}(u)$ is increasing in $u,$ and $u\frac{d}{du}\left(ur_X(u)\right)$ or $u\frac{d}{du}\left(ur_Y(u)\right)$is decreasing in $u,$ then, $\mbox{\boldmath $\frac{1}{\sigma}$}\stackrel{p}{\preceq} \mbox{\boldmath $\frac{1}{\xi}$}$ $\Rightarrow$ $X_{1:n}\leq_{hr}Y_{1:n}$,
\item[iv)] either $ur_{X}(u)$ or $ur_{Y}(u)$ is increasing in $u,$ and $u^2\frac{d}{du}\left(r_{X}(u)\right)$ or $u^2\frac{d}{du}\left(r_{Y}(u)\right)$is decreasing in $u,$ then, $\mbox{\boldmath $\frac{1}{\sigma}$}\stackrel{rm}{\preceq} \mbox{\boldmath $\frac{1}{\xi}$}$ $\Rightarrow$ $X_{1:n}\leq_{hr}Y_{1:n}.$
\end{enumerate}
\end{t1}
{\bf Proof:} If $X\leq_{hr}Y$, then $r_{X}\left(\frac{t-\lambda_k}{\sigma_k}\right)\ge r_{Y}\left(\frac{t-\lambda_k}{\sigma_k}\right)$, which results in 
$$\sum_{k=1}^n\frac{r_{X}\left(\frac{t-\lambda_k}{\sigma_k}\right)}{\sigma_k}\ge \sum_{k=1}^n\frac{r_{Y}\left(\frac{t-\lambda_k}{\sigma_k}\right)}{\sigma_k}.$$ 
So, to prove the result it suffices to prove that $\sum_{k=1}^n\frac{r_{Y}\left(\frac{t-\lambda_k}{\sigma_k}\right)}{\sigma_k}\ge \sum_{k=1}^n\frac{r_{Y}\left(\frac{t-\lambda_k}{\xi_k}\right)}{\xi_k}.$\\
$Proof~of~i)$:  Here we have to prove that $\sum_{k=1}^n\frac{r_{Y}\left(\frac{t-\lambda_k}{\sigma_k}\right)}{\sigma_k}=\sum_{k=1}^n p_{k} r_{Y}\left(p_{k}(t-\lambda_k)\right)=\Psi(\mbox{\boldmath $p$})$, is s-concave in $\mbox{\boldmath $p$}$. Now, if $\mbox{\boldmath $\lambda$},\mbox{\boldmath $\sigma$}, \mbox{\boldmath $\xi$}\in \mathcal{D}_+ \left(or~\mathcal{E}_+\right)$ then $p_i\left(t-\lambda_i\right)\leq (\ge) p_j\left(t-\lambda_j\right).$ Thus, the concavity of the function $ur_Y(u)$ gives
 $$\frac{d}{du}\left[ur_{Y}(u)\right]_{u=p_{i}(t-\lambda_i)}-\frac{d}{du}\left[ur_{Y}(u)\right]_{u=p_{j}(t-\lambda_j)}\ge (\leq) 0.$$
Now, differentiating $\Psi(\mbox{\boldmath $p$})$ with respect to $p_i$, we get
\begin{equation}\label{e22}
\frac{\partial \Psi}{\partial p_i}=\frac{d}{du}\left[ur_Y(u)\right]_{u=p_i(t-\lambda_i)},
\end{equation}
proving that $\frac{\partial \Psi}{\partial p_i}-\frac{\partial \Psi}{\partial p_j}\ge (\le) 0.$ Thus, by Lemma \ref{l7} (Lemma \ref{l6}) $\Psi(\mbox{\boldmath $p$})$ is s-concave in $\mbox{\boldmath $p$}$. This proves the result. \\  
$Proof~of~ii):$ If $u r_{Y}(u)$ is increasing (decreasing) in $u,$ then from (\ref{e22}) it can be written that that $\Psi(\mbox{\boldmath $p$})$ is increasing (decreasing) in $p_i.$ Again, as $ur_Y(u)$ is concave, by the previous result  $\Psi(\mbox{\boldmath $p$})$ is s-concave in $\mbox{\boldmath $p$}$. Thus, by Lemma \ref{l10}, the result is proved. \\
$Proof~of~iii):$ Let us assume that 
$$\sum_{k=1}^n\frac{ r_{Y}\left(\frac{t-\lambda_k}{\sigma_k}\right)}{\sigma_k}=\sum_{k=1}^n e^{p_k}r_{Y}\left(e^{p_k}(t-\lambda_k)\right)=\Psi(\mbox{\boldmath $p$}),$$ 
where $\mbox{\boldmath $p$}=\left(\log\frac{1}{\sigma_1},\log\frac{1}{\sigma_2},...,\log\frac{1}{\sigma_n}\right)$. So, by Lemma \ref{l8}, it is suffices to prove that $\Psi(\mbox{\boldmath $p$})$ is increasing in each $p_i$ and s-concave in $\mbox{\boldmath $p$}$. Now, differentiating $\Psi(\mbox{\boldmath $p$})$ with respect to $p_i$, we get,  
\begin{eqnarray*}
\frac{\partial \Psi}{\partial p_i}&=&e^{p_i}\left[r_{Y}\left(e^{p_i}(t-\lambda_i)\right)\right]+e^{p_i}(t-\lambda_i)r_{Y}\left(e^{p_i}(t-\lambda_i)\right)\\&=&e^{p_i}\frac{d}{du}\left[ur_{Y}(u)\right]_{u=e^{p_i}(t-\lambda_i)}\\&=&\frac{1}{t-\lambda_i}\left[u\frac{d}{du}\left(ur_Y(u)\right)\right]_{u=e^{p_i}(t-\lambda_i)}\\&>&0.
\end{eqnarray*} 
The last inequality follows from the fact that $ur_Y(u)$ is increasing in $u$. Again, as $u\frac{d}{du}\left(ur_Y(u)\right)$is decreasing in $u,$ and $\mbox{\boldmath $\lambda$},\mbox{\boldmath $\sigma$}\in \mathcal{D}_+ (or \mathcal{E}_+),$ then for all $i\leq j$, 
\begin{eqnarray*}
\frac{\partial \Psi}{\partial p_i}-\frac{\partial \Psi}{\partial p_j}&=&\frac{1}{t-\lambda_i}\left[u\frac{d}{du}\left(ur_Y(u)\right)\right]_{u=e^{p_i}(t-\lambda_i)}-\frac{1}{t-\lambda_j}\left[u\frac{d}{du}\left(ur_Y(u)\right)\right]_{u=e^{p_j}(t-\lambda_j)}\\&\ge (\le)&0.
\end{eqnarray*} Hence, the result is followed from Lemma \ref{l7} (Lemma \ref{l6}).\\ 
$Proof~of~iv)$ Let us assume that $\Psi(\mbox{\boldmath $\sigma$})=\sum_{k=1}^n\frac{r_{Y}\left(\frac{t-\lambda_k}{\sigma_k}\right)}{\sigma_k}.$ Now, differentiating $\Psi(\mbox{\boldmath $\sigma$})$ with respect to $\sigma_i$ and considering the fact that $ur_Y(u)$ is increasing in $u$, it can be written that, 
\begin{eqnarray*}
\frac{\partial \Psi}{\partial \sigma_i}&=&\frac{-1}{\sigma_{i}^{2}}r_{Y}\left(\frac{t-\lambda_i}{\sigma_i}\right)-\left(\frac{t-\lambda_i}{\sigma_i^3}\right) r^{'}_{Y}\left(\frac{t-\lambda_i}{\sigma_i}\right)\\&=&\frac{-1}{(t-\lambda_i)^2}\left(\frac{t-\lambda_i}{\sigma_i}\right)^2 \frac{d}{du}\left[ur_{Y}(u)\right]_{u=\frac{t-\lambda_i}{\sigma_i}}\\&\le& 0.
\end{eqnarray*} 
Again, as $\mbox{\boldmath $\lambda$},\mbox{\boldmath $\sigma$}\in \mathcal{D}_+ (or\; \mathcal{E}_+),$ and $u^2 \frac{d}{du}\left[ur_{Y}(u)\right]$ is decreasing in $u,$ it can be easily shown that $\frac{\partial \Psi}{\partial \sigma_i}-\frac{\partial \Psi}{\partial \sigma_j}\le (\ge) 0.$ Thus, by Lemma \ref{l6} (Lemma \ref{l7}), $\Psi(\mbox{\boldmath $\sigma$})$ is s-concave and hence the result is proved by Lemma \ref{l9}.   \hfill$\Box$
\begin{t1}\label{th4a}
Let $X_1,X_2,...,X_n$ be a set of independent random variables such that $X_i\sim $LS$\left(\lambda_i,\sigma_i,F\right)$, $i=1,2,...,n$. Let $Y_1,Y_2,...,Y_n$ be another set of independent random variables such that $Y_i\sim $LS$\left(\lambda_i,\xi_i,G\right),~i=1,2,...,n$. Suppose that $\mbox{\boldmath $\sigma$},\mbox{\boldmath $\xi$},\mbox{\boldmath $\lambda$}\in \mathcal{D}_+$ (or $\mathcal{E}_+$). Then, if $X\leq_{hr}Y$ and
\begin{enumerate}
\item[i)] either $ur_{X}(u)$ or $ur_{Y}(u)$ is convex in $u,$ then, $\mbox{\boldmath $\frac{1}{\sigma}$}\stackrel{m}{\succeq} \mbox{\boldmath $\frac{1}{\xi}$}$ $\Rightarrow$ $X_{1:n}\leq_{hr}Y_{1:n}$,
\item[ii)] either $ur_{X}(u)$ or $ur_{Y}(u)$ is increasing (decreasing) and convex in $u,$ then, $\mbox{\boldmath $\frac{1}{\sigma}$}\succeq_w(\stackrel{w}{\succeq}) \mbox{\boldmath $\frac{1}{\xi}$}$ $\Rightarrow$ $X_{1:n}\leq_{hr}Y_{1:n}$,
\item[iii)] either $ur_{X}(u)$ or $ur_{Y}(u)$ is decreasing in $u,$ and $u\frac{d}{du}\left(ur_X(u)\right)$ or $u\frac{d}{du}\left(ur_Y(u)\right)$is increasing in $u,$ then, $\mbox{\boldmath $\frac{1}{\sigma}$}\stackrel{p}{\succeq} \mbox{\boldmath $\frac{1}{\xi}$}$ $\Rightarrow$ $X_{1:n}\leq_{hr}Y_{1:n}$,
\item[iv)] either $ur_{X}(u)$ or $ur_{Y}(u)$ is decreasing in $u,$ and $u^2\frac{d}{du}\left(r_{X}(u)\right)$ or $u^2\frac{d}{du}\left(r_{Y}(u)\right)$is increasing in $u,$ then, $\mbox{\boldmath $\frac{1}{\sigma}$}\stackrel{rm}{\succeq} \mbox{\boldmath $\frac{1}{\xi}$}$ $\Rightarrow$ $X_{1:n}\leq_{hr}Y_{1:n}.$
\end{enumerate}
\end{t1}
{\bf Proof:} Proof of the theorem follows from the previous theorem and hence is omitted. \hfill$\Box$\\
\hspace*{0.3 in} In the next theorem, it is shown that under certain restrictions, there exists R-hr ordering between two minimum order statistics obtained from two different $LS$ families, when one set of scale parameters majorizes the other. It is to be mentioned here that, due to mathematical complexities we cannot proceed in full generality, where two $LS$ families could be generated from two different baseline distributions and hence we consider that both $LS$ families are generated from the same baseline distribution. 

%Theorem~\ref{th4}(i) guarantees that for series systems of components having independent $LS$ distributed lifetimes with common location parameter vector, the majorized scale parameter vector leads to larger systems lifetime in the sense of the hazard rate ordering. The following theorem shows that Theorem~\ref{th4} can be extended up to likelihood ratio ordering, although under more restrictive conditions.
\begin{t1}\label{th5}
Let $X_1,X_2,...,X_n$ be a set of independent random variables such that $X_i\sim $LS$\left(\lambda_i,\sigma_i,F\right)$, $i=1,2,...,n$. Let $Y_1,Y_2,...,Y_n$ be another set of independent random variables such that $Y_i\sim $LS$\left(\lambda_i,\xi_i,F\right)$, $i=1,2,...,n$. Suppose that $\mbox{\boldmath $\lambda$},\mbox{\boldmath $\sigma$},\mbox{\boldmath $\xi$}\in \mathcal{D}_+$ (or $\mathcal{E}_+$). If $r_{X}(u)$ is $r$-concave for $r=1,2$ and $3$ then, $\mbox{\boldmath $\frac{1}{\sigma}$}\stackrel{m}{\preceq} \mbox{\boldmath $\frac{1}{\xi}$}$ $\Rightarrow$ $X_{1:n}\le_{R-hr}Y_{1:n}$.
\end{t1}
{\bf Proof:} To prove the result we have only to prove that $g(t)=\frac{s_{1:n}(t)}{r_{1:n}(t)}$ is decreasing in $t$. Now, 
\begin{equation*}
g^{'}(t)=\frac{d}{dt}\left[\frac{\sum_{k=1}^n \frac{1}{\xi_k}r_{X}\left(\frac{t-\lambda_k}{\xi_k}\right)}{\sum_{k=1}^n\frac{1}{\sigma_k} r_{X}\left(\frac{t-\lambda_k}{\sigma_k}\right)}\right]\\\stackrel{sign}{=}\frac{\sum_{k=1}^n \frac{1}{\xi^{2}_{k}}r^{'}_{X}\left(\frac{t-\lambda_k}{\xi_k}\right)}{\sum_{k=1}^n \frac{1}{\xi_{k}}r^{'}_{X}\left(\frac{t-\lambda_k}{\xi_k}\right)}-\frac{\sum_{k=1}^n \frac{1}{\sigma^{2}_{k}}r^{'}_{X}\left(\frac{t-\lambda_k}{\sigma_k}\right)}{\sum_{k=1}^n \frac{1}{\sigma_{k}}r^{'}_{X}\left(\frac{t-\lambda_k}{\sigma_k}\right)}.
\end{equation*}
So, it is to be proved that
$$\frac{\sum_{k=1}^n \frac{r^{'}_{X}\left(\frac{t-\lambda_k}{\sigma_k}\right)}{\sigma^{2}_{k}}}{\sum_{k=1}^n \frac{r^{'}_{X}\left(\frac{t-\lambda_k}{\sigma_k}\right)}{\sigma_{k}}}=\frac{\sum_{k=1}^n p^{2}_{k}r^{'}_{X}\left(p_k(t-\lambda_k)\right)}{\sum_{k=1}^n p_k r_{X}\left(p_k(t-\lambda_k)\right)}=\frac{\sum_{k=1}^n u(p_k,\lambda_k,t) v(p_k,\lambda_k,t)}{\sum_{k=1}^n u(p_k,\lambda_k,t)}=\frac{\sum_{k=1}^n u_k v_k}{\sum_{k=1}^n u_k}=\Psi(\mbox{\boldmath $p$})~(say),$$ 
with $\mbox{\boldmath $p$}=\left(\frac{1}{\sigma_1},\frac{1}{\sigma_2},...,\frac{1}{\sigma_n}\right),$ is s-concave in $\mbox{\boldmath $p$},$where $u(p_k,\lambda_k,t)=p_k r_{X}\left(p_k(t-\lambda_k)\right)=u_k$ (say) and $v(p_k,\lambda_k,t)=p_k\frac{r^{'}_{X}\left(p_k(t-\lambda_k)\right)}{r_{X}\left(p_k(t-\lambda_k)\right)}=v_k$ (say).  Now, differentiating $\Psi(\mbox{\boldmath $p$})$ with respect to $p_i$ we get, 
\begin{equation}\label{e23}
\frac{\partial \Psi}{\partial p_i}\stackrel{sign}{=}\sum_{k=1}^n u_k\left(v_i\frac{\partial u_i}{\partial p_i}+u_i \frac{\partial v_i}{\partial p_i}\right)-\sum_{k=1}^n u_kv_k\left(\frac{\partial u_i}{\partial p_i}\right).
\end{equation}  
 Now, differentiating $u_k$ and $v_k$ with respect to $p_k$, it can be written that \\
$\frac{\partial u_k}{\partial p_k}=u_k w_k;~w_k=\frac{1+v_k (t-\lambda_k)}{p_k}$, and\\
 $\frac{\partial v_k}{\partial p_k}=v_k\left[\frac{1}{p_k}+(t-\lambda_k)\frac{r^{''}_{X}\left(p_k(t-\lambda_k)\right)}{r^{'}_{X}\left(p_k(t-\lambda_k)\right)}-(t-\lambda_k)\frac{r^{'}_{X}\left(p_k(t-\lambda_k)\right)}{r_{Y}\left(p_k(t-\lambda_k)\right)}\right]=v_k s_k$ (say),\\
where $s_k=\frac{1}{p_k}+(t-\lambda_k)\frac{r^{''}_{X}\left(p_k(t-\lambda_k)\right)}{r^{'}_{X}\left(p_k(t-\lambda_k)\right)}-(t-\lambda_k)\frac{r^{'}_{X}\left(p_k(t-\lambda_k)\right)}{r_{Y}\left(p_k(t-\lambda_k)\right)}.$\\
Using (\ref{e23}) and expressions for $\frac{\partial u_k}{\partial p_k}$ and $\frac{\partial v_k}{\partial p_k}$, it can be shown that, 
$$ \frac{\partial \Psi}{\partial p_i}\stackrel{sign}{=}\left(\sum_{k=1}^n u_k\right)\left(u_i v_i w_i+u_i v_i s_i\right)-\left(\sum_{k=1}^n u_k\right) v_k u_i w_i,$$ which gives, for all $i\leq j$
\begin{equation}\label{e24}
\frac{\partial \Psi}{\partial p_i}-\frac{\partial \Psi}{\partial p_j}\stackrel{sign}{=}\left(\sum_{k=1}^n u_k\right)\left[u_i v_i(w_i+s_i)-u_j v_j(w_j+s_j)\right]-\left(\sum_{k=1}^n u_k v_k\right)\left(u_i w_i-u_j w_j\right),
\end{equation}
with 
\begin{equation}\label{e25}
\begin{split}
u_i v_i(w_i+s_i)-u_j v_j(w_j+s_j)&=2\left[p_i r^{'}_{X}\left(p_i(t-\lambda_i)\right)-p_j r^{'}_{X}\left(p_j(t-\lambda_j)\right)\right]\\&\quad +\left[p^{2}_{i}(t-\lambda_i) r^{''}_{X}\left(p_i(t-\lambda_i)\right)-p^{2}_{j} r^{''}_{X}\left(p_j(t-\lambda_j)\right)\right].
\end{split}
\end{equation}
Now as $\mbox{\boldmath $\lambda$},\mbox{\boldmath $\sigma$}\in \mathcal{D}_+, (or~\mathcal{E_+})$, then for all $i\leq j$, $p_i,<p_j$, which gives $p_i(t-\lambda_i)\le (\ge) p_j(t-\lambda_j)$ and $p^{2}_{i}(t-\lambda_i)\le (\ge) p^{2}_{j}(t-\lambda_j).$ Moreover, as $r_{X}(u)$ is s-concave for $s=1,2$, it is decreasing and concave and hence it can be concluded that $p_i r^{'}_{X}\left(p_i(t-\lambda_i)\right)\ge (\le)p_j r^{'}_{X}\left(p_j(t-\lambda_j)\right)$. So, the first term of (\ref{e25}) is positive (negative).\\
Again, as $r_{X}(u)$ is also $3$-concave and thus $r^{''}_{X} (u)$ decreasing in $u$, then for all $i\leq j$ it can be written that, $r^{''}_{X} (p_i(t-\lambda_i))\ge (\le)r^{''}_{X} (p_j(t-\lambda_j))$, which gives $p^{2}_{i}(t-\lambda_i) r^{''}_{X}\left(p_i(t-\lambda_i)\right)\ge (\le)p^{2}_{j}(t-\lambda_j)r^{''}_{X}\left(p_j(t-\lambda_j)\right)$. So the second term of (\ref{e25}) is also positive (negative). Thus the first term of (\ref{e24}) is positive (negative). \\
Now, following the arguments as before it can be easily shown that,
$$u_i w_i-u_j w_j=\left[r_{X}\left(p_i(t-\lambda_i)\right)-r_{X}\left(p_j(t-\lambda_j)\right)\right]+\left[p_i(t-\lambda_i)r^{'}_{X}\left(p_i(t-\lambda_i)\right)-p_j(t-\lambda_j)r^{'}_{X}\left(p_j(t-\lambda_j)\right)\right]\ge (\le)0.$$
Now as $r_X(u)$ is decreasing in $u$, $v_k=p_k\frac{r^{'}_{X}\left(p_k(t-\lambda_k)\right)}{r_{X}\left(p_k(t-\lambda_k)\right)}$ is negative. So, the second term of (\ref{e24}) is also positive (negative). Hence, $\frac{\partial \Psi}{\partial p_i}-\frac{\partial \Psi}{\partial p_j}\ge (\le) 0$. Thus by Lemma \ref{l7} (Lemma \ref{l6}) it can be concluded that $\Psi(\mbox{\boldmath $p$})$ is s-concave in $\mbox{\boldmath $p$}$. This proves the result.\hfill$\Box$\\\\
\hspace*{0.3 in} Systems are also compared stochastically in situations where the components are from multiple-outlier models. This is to be mentioned here that, a multiple-outlier model is a set of independent random variables $X_1,..., X_n$ of which $X_i\stackrel{st}{=}X, i=1,...,n_1$ and $X_i\stackrel{st}{=}Y, i=n_1+1,...,n$ where $1\leq n_1<n$ and $X_i\stackrel{st}{=}X$ means that cdf of $X_i$ is same as that of $X.$ In other words, the set of independent random variables $X_1,..., X_n$ is said to constitute a multiple-outlier model if two sets of random variables $\left(X_1, X_2,\ldots, X_{n_1}\right)$ and $\left(X_{n_1+1}, X_{n_1+2},\ldots, X_{n_1+n_2}\right)$ (where $n_1+n_2=n$), are homogenous among themselves and heterogenous between themselves. For more details on multiple-outlier models, readers may refer to  Kochar and Xu~(\cite{ko1}), Zhao and Balakrishnan~(\cite{zb11}), Balakrishnan and Torrado~(\cite{ba}), Zhao and Zhang~\cite{zh1},  Kundu \emph{et al.}~\cite{kun1}, Kundu and Chowdhury~\cite{kun2} and the references there in. Two theorems given below show that, under certain conditions, majorized scale parameter vectors of $LS$ distributed components of series system for multiple outlier model, leads to smaller system life in terms of R-hr ordering.
\begin{t1}\label{th9}
Let $\left\{X_1, X_2, \ldots, X_n\right\}$ and $\left\{Y_1, Y_2, \ldots, Y_n\right\}$ be two sets of independent random variables each following the multiple outlier $LS$ model such that $X_i\sim $LS$\left(\lambda_1,\sigma_1,F\right)$ and $Y_i\sim $LS$\left(\lambda_1,\xi_1,F\right)$ for $i=1,2,\ldots,n_1$ and $X_i\sim $LS$\left(\lambda_2,\sigma_2,F\right)$ and $Y_i\sim $LS$\left(\lambda_2,\xi_2,F\right)$ for $i=n_1+1,n_1+2,\ldots,n_1+n_2=n$. Suppose that $\left(\lambda_1, \lambda_2\right), \left(\sigma_1, \sigma_2\right), \left(\xi_1, \xi_2\right)\in \mathcal{D}_+ \;(or\; \mathcal{E}_+)$. Now if $ur_X(u)$ is decreasing in $u$ and $r_X(u)$ is log-concave and $2$-log-convex in $u$, then 
$$(\underbrace{\frac{1}{\sigma_1},\frac{1}{\sigma_1},\ldots,\frac{1}{\sigma_1},}_{n_1} \underbrace{\frac{1}{\sigma_2},\frac{1}{\sigma_2},\ldots,\frac{1}{\sigma_2}}_{n_2})\stackrel{m}{\preceq}
(\underbrace{\frac{1}{\xi_1},\frac{1}{\xi_1},\ldots,\frac{1}{\xi_1},}_{n_1} \underbrace{\frac{1}{\xi_2},\frac{1}{\xi_2},\ldots,\frac{1}{\xi_2}}_{n_2})\Rightarrow X_{1:n}\geq_{R-hr}Y_{1:n}.$$ 
\end{t1}
{\bf Proof:} To prove the result we have only to prove that $g(t)=\frac{s_{1:n}(t)}{r_{1:n}(t)}$ is increasing in $t$. Now, 
\begin{equation*}
g^{'}(t)=\frac{d}{dt}\left[\frac{\sum_{k=1}^n \frac{1}{\xi_k}r_{X}\left(\frac{t-\lambda_k}{\xi_k}\right)}{\sum_{k=1}^n\frac{1}{\sigma_k} r_{X}\left(\frac{t-\lambda_k}{\sigma_k}\right)}\right]\\\stackrel{sign}{=}\frac{\sum_{k=1}^n \frac{1}{\xi^{2}_{k}}r^{'}_{X}\left(\frac{t-\lambda_k}{\xi_k}\right)}{\sum_{k=1}^n \frac{1}{\xi_{k}}r^{'}_{X}\left(\frac{t-\lambda_k}{\xi_k}\right)}-\frac{\sum_{k=1}^n \frac{1}{\sigma^{2}_{k}}r^{'}_{X}\left(\frac{t-\lambda_k}{\sigma_k}\right)}{\sum_{k=1}^n \frac{1}{\sigma_{k}}r^{'}_{X}\left(\frac{t-\lambda_k}{\sigma_k}\right)}.
\end{equation*}
So, to prove the result it is to be shown that
$$\frac{\sum_{k=1}^n \frac{r^{'}_{X}\left(\frac{t-\lambda_k}{\sigma_k}\right)}{\sigma^{2}_{k}}}{\sum_{k=1}^n \frac{r^{'}_{X}\left(\frac{t-\lambda_k}{\sigma_k}\right)}{\sigma_{k}}}=\frac{\sum_{k=1}^n p^{2}_{k}r^{'}_{X}\left(p_k(t-\lambda_k)\right)}{\sum_{k=1}^n p_k r_{X}\left(p_k(t-\lambda_k)\right)}=\frac{\sum_{k=1}^n u(p_k,\lambda_k,t) v(p_k,\lambda_k,t)}{\sum_{k=1}^n u(p_k,\lambda_k,t)}=\frac{\sum_{k=1}^n u_k v_k}{\sum_{k=1}^n u_k}=\Psi(\mbox{\boldmath $p$})~(say),$$ 
with $\mbox{\boldmath $p$} =\left(\underbrace{p_1,\ldots,p_1}_{n_1},\underbrace{p_2,\ldots,p_2}_{n_2}\right)$ $=\left(\underbrace{\frac{1}{\sigma_1},\frac{1}{\sigma_1},\ldots,\frac{1}{\sigma_1}}_{n_1},\underbrace{\frac{1}{\sigma_2}, \frac{1}{\sigma_2}\ldots\frac{1}{\sigma_2}}_{n_2}\right),$ is s-convex in $\mbox{\boldmath $p$},$ where $u(p_k,\lambda_k,t)=p_k r_{X}\left(p_k(t-\lambda_k)\right)=u_k$ (say) and $v(p_k,\lambda_k,t)=p_k\frac{r^{'}_{X}\left(p_k(t-\lambda_k)\right)}{r_{X}\left(p_k(t-\lambda_k)\right)}=v_k$ (say).\\
Let $i\leq j$. Now three cases may arise:\\
$Case (i):$ If $1\leq i,j\leq n_1$, $i.e.$, if $p_i=p_j=p_1$ and $\lambda_i=\lambda_j=\lambda_1$, then 
$$\frac{\partial \Psi}{\partial p_i}-\frac{\partial \Psi}{\partial p_j}=\frac{\partial \Psi}{\partial p_1}-\frac{\partial \Psi}{\partial p_1}=0.$$
$Case (ii):$ If $n_1+1\leq i,j\leq n$, $i.e.$, if $p_i=p_j=p_2$ and $\lambda_i=\lambda_j=\lambda_2$, then 
$$\frac{\partial \Psi}{\partial p_i}-\frac{\partial \Psi}{\partial p_j}=\frac{\partial \Psi}{\partial p_2}-\frac{\partial \Psi}{\partial p_2}=0.$$
$Case (iii):$ If $1\leq i\leq n_1$ and $n_1+1\leq j\leq n$, then $p_i=p_1$, $p_j=p_2$, $\lambda_i=\lambda_1$ and $\lambda_j=\lambda_2$. Then, 
\begin{equation*}
\frac{\partial \Psi}{\partial p_i}\left(n_1u_1+n_2u_2\right)^2=n_2u_2\frac{\partial u_1}{\partial p_1}\left(v_1-v_2\right)+w_1\left(n_1u_1+n_2u_2\right),
\end{equation*}
and
\begin{equation*}
\frac{\partial \Psi}{\partial p_j}\left(n_1u_1+n_2u_2\right)^2=n_1u_1\frac{\partial u_2}{\partial p_2}\left(v_2-v_1\right)+w_2\left(n_1u_1+n_2u_2\right),
\end{equation*}
where, for $i=1,2$, $w_i=u_i\frac{\partial v_i}{\partial p_i}$. Thus 
\begin{equation}\label{1}
\frac{\partial \Psi}{\partial p_i}-\frac{\partial \Psi}{\partial p_j}\stackrel{sign}{=} \left(v_1-v_2\right)\left(n_2u_2\frac{\partial u_1}{\partial p_1}+n_1u_1\frac{\partial u_2}{\partial p_2}\right)+\left(w_1-w_2\right)\left(n_1u_1+n_2u_2\right).
\end{equation}
Now, as $ur_X(u)$ is decreasing in $u$, then for $i=1,2$, $\frac{\partial u_i}{\partial p_i}=\frac{d}{du}\left(ur_X(u)\right)|_{u=\left(t-\lambda_i\right)p_i}<0$. Now, if $\left(\lambda_1, \lambda_2\right), \left(\sigma_1, \sigma_2\right)\in \mathcal{D}_+ \;(or \mathcal{E}_+)$, then $p_1\le(\ge)p_2$ and $p_1\left(t-\lambda_1\right)\le (\ge)p_2\left(t-\lambda_2\right)$. Thus, as $r_X(u)$ is log-concave in $u$, by considering the fact that $r_X(u)$ is decreasing in $u$, it can be written that
$$-p_1\frac{r_X^{'}\left(p_1\left(t-\lambda_1\right)\right)}{r_X\left(p_1\left(t-\lambda_1\right)\right)}\le(\ge)-p_2\frac{r_X^{'}\left(p_2\left(t-\lambda_2\right)\right)}{r_X\left(p_2\left(t-\lambda_2\right)\right)},$$
giving that $v_1\ge (\le) v_2$. So, the first term of (\ref{1}) is negative (positive). Again, for $i=1,2$, 
\begin{eqnarray*}
w_i&=&p_ir_X\left(t-\lambda_i\right)\frac{\partial}{\partial p_i}\left(\frac{r_X^{'}\left(p_i\left(t-\lambda_i\right)\right)}{r_X\left(p_i\left(t-\lambda_i\right)\right)}\right)\\
&=&p_i\left(t-\lambda_i\right)r_X\left(t-\lambda_i\right)\frac{\partial^2}{\partial u^2}\left(\log r_X\left(u\right)\right)|_{u=p_i\left(t-\lambda_i\right)}.
\end{eqnarray*}
So, considering the fact that $r_X(u)$ is log-concave giving $\frac{\partial^2}{\partial u^2}\log r_X\left(u\right)<0$, and $r_X(u)$ is also $2$-log-convex, it can be written that, 
$$-p_1\left(t-\lambda_1\right)r_X\left(t-\lambda_1\right)\frac{\partial^2}{\partial u^2}\left(\log r_X\left(u\right)\right)|_{u=p_1\left(t-\lambda_1\right)}\geq (\leq) -p_2\left(t-\lambda_2\right)r_X\left(t-\lambda_2\right)\frac{\partial^2}{\partial u^2}\left(\log r_X\left(u\right)\right)|_{u=p_2\left(t-\lambda_2\right)},$$ yielding that $w_1\leq (\geq)w_2$. So, the second term of (\ref{1}) is also negative (positive). Thus, as $\frac{\partial \Psi}{\partial p_i}-\frac{\partial \Psi}{\partial p_j}\leq(\geq) 0$, by Lemma \ref{l7} (Lemma \ref{l6}) it can be concluded that $\Psi(\mbox{\boldmath $p$})$ is s-convex in $\mbox{\boldmath $p$}$. This proves the result.\hfill$\Box$\\
The next theorem can be proved in the similar line as of previous theorem.
\begin{t1}\label{th10}
Let $\left\{X_1, X_2, \ldots, X_n\right\}$ and $\left\{Y_1, Y_2, \ldots, Y_n\right\}$ be two sets of independent random variables each following the multiple outlier $LS$ model such that $X_i\sim $LS$\left(\lambda_1,\sigma_1,F\right)$ and $Y_i\sim $LS$\left(\lambda_1,\xi_1,F\right)$ for $i=1,2,\ldots,n_1$ and $X_i\sim $LS$\left(\lambda_2,\sigma_2,F\right)$ and $Y_i\sim $LS$\left(\lambda_2,\xi_2,F\right)$ for $i=n_1+1,n_1+2,\ldots,n_1+n_2=n$. Suppose that $\left(\lambda_1, \lambda_2\right), \left(\sigma_1, \sigma_2\right), \left(\xi_1, \xi_2\right)\in \mathcal{D}_+ \;(or\; \mathcal{E}_+)$. Now if $ur_X(u)$ is increasing in $u$ and $r_X(u)$ is decreasing, log-concave and $2$-log-concave in $u$, then 
$$(\underbrace{\frac{1}{\sigma_1},\frac{1}{\sigma_1},\ldots,\frac{1}{\sigma_1},}_{n_1} \underbrace{\frac{1}{\sigma_2},\frac{1}{\sigma_2},\ldots,\frac{1}{\sigma_2}}_{n_2})\stackrel{m}{\succeq}
(\underbrace{\frac{1}{\xi_1},\frac{1}{\xi_1},\ldots,\frac{1}{\xi_1},}_{n_1} \underbrace{\frac{1}{\xi_2},\frac{1}{\xi_2},\ldots,\frac{1}{\xi_2}}_{n_2})\Rightarrow X_{1:n}\geq_{R-hr}Y_{1:n}.$$
\end{t1}
Now, combining Theorem \ref{th4} i) and Theorem \ref{th9} it can be shown that, for multiple outlier $LS$ model, majorized scale parameter vector of minimum order statistics leads to smaller system life in terms of lr ordering. The statement of the theorem is given below.
\begin{t1}\label{th11}
Let $\left\{X_1, X_2, \ldots, X_n\right\}$ and $\left\{Y_1, Y_2, \ldots, Y_n\right\}$ be two sets of independent random variables each following the multiple outlier $LS$ model such that $X_i\sim $LS$\left(\lambda_1,\sigma_1,F\right)$ and $Y_i\sim $LS$\left(\lambda_1,\xi_1,F\right)$ for $i=1,2,\ldots,n_1$ and $X_i\sim $LS$\left(\lambda_2,\sigma_2,F\right)$ and $Y_i\sim $LS$\left(\lambda_2,\xi_2,F\right)$ for $i=n_1+1,n_1+2,\ldots,n_1+n_2=n$. Suppose that $\left(\lambda_1, \lambda_2\right), \left(\sigma_1, \sigma_2\right), \left(\xi_1, \xi_2\right)\in \mathcal{D}_+ \;(or\; \mathcal{E}_+)$. Now if $ur_X(u)$ is decreasing and concave in $u$ and $r_X(u)$ is log-concave and $2$-log-convex in $u$, then 
$$(\underbrace{\frac{1}{\sigma_1},\frac{1}{\sigma_1},\ldots,\frac{1}{\sigma_1},}_{n_1} \underbrace{\frac{1}{\sigma_2},\frac{1}{\sigma_2},\ldots,\frac{1}{\sigma_2}}_{n_2})\stackrel{m}{\preceq}
(\underbrace{\frac{1}{\xi_1},\frac{1}{\xi_1},\ldots,\frac{1}{\xi_1},}_{n_1} \underbrace{\frac{1}{\xi_2},\frac{1}{\xi_2},\ldots,\frac{1}{\xi_2}}_{n_2})\Rightarrow X_{1:n}\leq_{lr}Y_{1:n}.$$ 
\end{t1}
Again combining Theorem \ref{th4a} i) and Theorem \ref{th10} the following theorem can be obtained.
\begin{t1}\label{th12}
Let $\left\{X_1, X_2, \ldots, X_n\right\}$ and $\left\{Y_1, Y_2, \ldots, Y_n\right\}$ be two sets of independent random variables each following the multiple outlier $LS$ model such that $X_i\sim $LS$\left(\lambda_1,\sigma_1,F\right)$ and $Y_i\sim $LS$\left(\lambda_1,\xi_1,F\right)$ for $i=1,2,\ldots,n_1$ and $X_i\sim $LS$\left(\lambda_2,\sigma_2,F\right)$ and $Y_i\sim $LS$\left(\lambda_2,\xi_2,F\right)$ for $i=n_1+1,n_1+2,\ldots,n_1+n_2=n$. Suppose that $\left(\lambda_1, \lambda_2\right), \left(\sigma_1, \sigma_2\right), \left(\xi_1, \xi_2\right)\in \mathcal{D}_+ \;(or\; \mathcal{E}_+)$. Now if $ur_X(u)$ is increasing and convex in $u$ and $r_X(u)$ is decreasing, log-concave and $2$-log-convex in $u$, then 
$$(\underbrace{\frac{1}{\sigma_1},\frac{1}{\sigma_1},\ldots,\frac{1}{\sigma_1},}_{n_1} \underbrace{\frac{1}{\sigma_2},\frac{1}{\sigma_2},\ldots,\frac{1}{\sigma_2}}_{n_2})\stackrel{m}{\succeq}
(\underbrace{\frac{1}{\xi_1},\frac{1}{\xi_1},\ldots,\frac{1}{\xi_1},}_{n_1} \underbrace{\frac{1}{\xi_2},\frac{1}{\xi_2},\ldots,\frac{1}{\xi_2}}_{n_2})\Rightarrow X_{1:n}\leq_{lr}Y_{1:n}.$$ 
\end{t1}
\hspace*{0.3 in} Now the question arises$~-~$what will happen if location parameter vector of one $LS$ family majorizes the other when the
scale parameter vector remains same? The next few theorems deal with such cases.\\
The next theorem shows that under some restrictions, majorization ordering between location parameter vectors of minimum order statistics from two different $LS$ family of distributions implies hr ordering between them.
\begin{t1}\label{th13}
Let $\left\{X_1, X_2, \ldots, X_n\right\}$ and $\left\{Y_1, Y_2, \ldots, Y_n\right\}$ be two sets of independent random variables each following $LS$ model such that $X_i\sim $LS$\left(\lambda_i,\sigma_i,F\right)$ and $Y_i\sim $LS$\left(\mu_i,\sigma_i,G\right)$ for $i=1,2,\ldots,n$. Suppose further that $\mbox{\boldmath $\lambda$},\mbox{\boldmath $\mu$},\mbox{\boldmath $\sigma$}\in \mathcal{D}_+$ (or $\mathcal{E}_+$). Then, if $X\leq_{hr}Y$, either $r_X(u)$ or $r_Y(u)$ is increasing (decreasing) in $u$ and either $u^2r_X^{'}(u)$ or $u^2r_Y^{'}(u)$ is decreasing (increasing) in $u$ then, $\mbox{\boldmath $\lambda$}\preceq_w (\succeq_w) \mbox{\boldmath $\mu$}$ $\Rightarrow$ $X_{1:n}\leq_{hr}Y_{1:n}$.
\end{t1}
{\bf Proof:} If $X\leq_{hr}Y$, then $r_{X}\left(\frac{t-\lambda_k}{\sigma_k}\right)\ge r_{Y}\left(\frac{t-\lambda_k}{\sigma_k}\right)$, which results in 
$$\sum_{k=1}^n\frac{r_{X}\left(\frac{t-\lambda_k}{\sigma_k}\right)}{\sigma_k}\ge \sum_{k=1}^n\frac{r_{Y}\left(\frac{t-\lambda_k}{\sigma_k}\right)}{\sigma_k}.$$ 
So, to prove the result it suffices to prove that $\sum_{k=1}^n\frac{r_{Y}\left(\frac{t-\lambda_k}{\sigma_k}\right)}{\sigma_k}\ge \sum_{k=1}^n\frac{r_{Y}\left(\frac{t-\mu_k}{\sigma_k}\right)}{\sigma_k}$, or equivalently it is required to show that $\sum_{k=1}^n\frac{r_{Y}\left(\frac{t-\mu_k}{\sigma_k}\right)}{\sigma_k}=\Psi(\mbox{\boldmath $\mu$})$ is decreasing (increasing) in each $\mu_i$ and s-concave (s-convex) in $\mbox{\boldmath $\mu$}$. Now, differentiating $\Psi$ with respect to $\mu_i$ it can be written that $\frac{\partial\Psi}{\partial\mu_i}=-\frac{1}{\sigma_i^2}r_Y^{'}\left(\frac{t-\mu_i}{\sigma_i}\right)$, which is clearly decreasing (increasing) in $\mu_i$ if $r_Y(u)$ is increasing (decreasing) in $u$. Again, if $\mbox{\boldmath $\mu$},\mbox{\boldmath $\sigma$}\in \mathcal{D}_+$, then for all $i\leq j$, $\mu_i\geq\mu_j$, $\sigma_i\geq\sigma_j$ giving $\frac{t-\mu_i}{\sigma_i}\leq\frac{t-\mu_j}{\sigma_j}$. So, if $r_Y(u)$ is increasing (decreasing)  and $u^2r_Y^{'}(u)$ is decreasing (increasing) in $u$, then
$$\frac{1}{(t-\mu_i)^2}\left(\frac{t-\mu_i}{\sigma_i}\right)^2r_Y^{'}\left(\frac{t-\mu_i}{\sigma_i}\right)\geq (\leq)\frac{1}{(t-\mu_j)^2}\left(\frac{t-\mu_j}{\sigma_j}\right)^2r_Y^{'}\left(\frac{t-\mu_j}{\sigma_j}\right)$$
giving that
$$\frac{\partial\Psi}{\partial\mu_i}-\frac{\partial\Psi}{\partial\mu_j}=\frac{1}{\sigma_j^2}r_Y^{'}\left(\frac{t-\mu_j}{\sigma_j}\right)-\frac{1}{\sigma_i^2}r_Y^{'}\left(\frac{t-\mu_i}{\sigma_i}\right)\leq (\geq) 0.$$
So, by Lemma \ref{l6} it can be concluded that $\Psi(\mbox{\boldmath $\mu$})$ is s-concave (s-convex) in $\mbox{\boldmath $\mu$}$. This proves the result.\\
When $\mbox{\boldmath $\mu$},\mbox{\boldmath $\sigma$}\in \mathcal{E}_+$, then the theorem can be proved in similar line. \hfill$\Box$\\
The theorem given below shows that for two $LS$ families of distribution functions generated from the same baseline distribution, there exists R-hr ordering between their corresponding minimum order statistics.
\begin{t1}\label{th14}
Let $\left\{X_1, X_2, \ldots, X_n\right\}$ and $\left\{Y_1, Y_2, \ldots, Y_n\right\}$ be two sets of independent random variables each following $LS$ model such that $X_i\sim $LS$\left(\lambda_i,\sigma_i,F\right)$ and $Y_i\sim $LS$\left(\mu_i,\sigma_i,F\right)$ for $i=1,2,\ldots,n$. Suppose further that $\mbox{\boldmath $\lambda$},\mbox{\boldmath $\mu$},\mbox{\boldmath $\sigma$}\in \mathcal{D}_+$ (or $\mathcal{E}_+$). Then, if $r_X(u)$ is increasing and concave (decreasing and convex) in $u$; $u^2r_X^{'}(u)$ is decreasing (increasing) in $u$ and $u^3r_X^{''}(u)$ is increasing (decreasing) in $u$ then, $\mbox{\boldmath $\lambda$}\preceq_w (\succeq_w) \mbox{\boldmath $\mu$}$ $\Rightarrow$ $X_{1:n}\geq_{R-hr}Y_{1:n}$.
\end{t1}
{\bf Proof:} To prove the result we have only to prove that $g(t)=\frac{s_{1:n}(t)}{r_{1:n}(t)}$ is increasing in $t$. Now, 
\begin{equation*}
g^{'}(t)=\frac{d}{dt}\left[\frac{\sum_{k=1}^n \frac{1}{\sigma_k}r_{X}\left(\frac{t-\mu_k}{\sigma_k}\right)}{\sum_{k=1}^n\frac{1}{\sigma_k} r_{X}\left(\frac{t-\lambda_k}{\sigma_k}\right)}\right]\\\stackrel{sign}{=}\frac{\sum_{k=1}^n \frac{1}{\sigma^{2}_{k}}r^{'}_{X}\left(\frac{t-\mu_k}{\sigma_k}\right)}{\sum_{k=1}^n \frac{1}{\xi_{k}}r^{'}_{X}\left(\frac{t-\mu_k}{\sigma_k}\right)}-\frac{\sum_{k=1}^n \frac{1}{\sigma^{2}_{k}}r^{'}_{X}\left(\frac{t-\lambda_k}{\sigma_k}\right)}{\sum_{k=1}^n \frac{1}{\sigma_{k}}r^{'}_{X}\left(\frac{t-\lambda_k}{\sigma_k}\right)}.
\end{equation*}
So, it is to be proved that
$$\frac{\sum_{k=1}^n \frac{r^{'}_{X}\left(\frac{t-\mu_k}{\sigma_k}\right)}{\sigma^{2}_{k}}}{\sum_{k=1}^n \frac{r^{'}_{X}\left(\frac{t-\mu_k}{\sigma_k}\right)}{\sigma_{k}}}=\Psi(\mbox{\boldmath $\mu$})$$
is increasing (decreasing) in each $\mu_i$ and s-convex (s-concave) in $\mbox{\boldmath $\mu$}$. Now differentiating $\Psi$ with respect to $\mu_i$, 
$$\frac{\partial\Psi}{\partial\mu_i}\left[\sum_{k=1}^n \frac{1}{\sigma_k}r_{X}\left(\frac{t-\mu_k}{\sigma_k}\right)\right]^2=-\left[\sum_{k=1}^n \frac{1}{\sigma_k}r_{X}\left(\frac{t-\mu_k}{\sigma_k}\right)\right]\frac{1}{\sigma_i^3}r_X^{''}\left(\frac{t-\mu_i}{\sigma_i}\right)+\left[\sum_{k=1}^n \frac{1}{\sigma_k^2}r_{X}^{'}\left(\frac{t-\mu_k}{\sigma_k}\right)\right]\frac{1}{\sigma_i^2}r_X^{'}\left(\frac{t-\mu_i}{\sigma_i}\right)$$
which is increasing (decreasing) in each $\mu_i$ if $r_X(u)$ is increasing and concave (decreasing and convex) in $u$. Again, 
\begin{equation}\label{2}
\begin{split}
\frac{\partial\Psi}{\partial\mu_i}-\frac{\partial\Psi}{\partial\mu_j}\stackrel{sign}{=}\left[\sum_{k=1}^n \frac{1}{\sigma_k}r_{X}\left(\frac{t-\mu_k}{\sigma_k}\right)\right]\left[\frac{1}{\sigma_j^3}r_X^{''}\left(\frac{t-\mu_j}{\sigma_j}\right)-\frac{1}{\sigma_i^3}r_X^{''}\left(\frac{t-\mu_i}{\sigma_i}\right)\right]&\\\quad+\left[\sum_{k=1}^n \frac{1}{\sigma_k^2}r_{X}^{'}\left(\frac{t-\mu_k}{\sigma_k}\right)\right]\left[\frac{1}{\sigma_i^2}r_X^{'}\left(\frac{t-\mu_i}{\sigma_i}\right)-\frac{1}{\sigma_j^2}r_X^{'}\left(\frac{t-\mu_j}{\sigma_j}\right)\right].
\end{split}
\end{equation}
Now if $\mbox{\boldmath $\mu$},\mbox{\boldmath $\sigma$}\in \mathcal{D}_+$, then for all $i\leq j$, $t-\mu_i\leq t-\mu_j$, $\sigma_i\geq\sigma_j$ and $\frac{t-\mu_i}{\sigma_i}\leq\frac{t-\mu_j}{\sigma_j}$. So, as $u^2r_X^{'}(u)$ is decreasing (increasing) in $u$ then for all $i\leq j$,
$$\left(\frac{t-\mu_i}{\sigma_i}\right)^2r_X^{'}\left(\frac{t-\mu_i}{\sigma_i}\right)\geq (\leq) \left(\frac{t-\mu_j}{\sigma_j}\right)^2r_X^{'}\left(\frac{t-\mu_j}{\sigma_j}\right).$$
As, $r_X(u)$ is increasing (decreasing) in $u$, then from the above relation it can be written that
$$\frac{1}{\left(t-\mu_i\right)^2}\left(\frac{t-\mu_i}{\sigma_i}\right)^2r_X^{'}\left(\frac{t-\mu_i}{\sigma_i}\right)\geq (\leq)\frac{1}{\left(t-\mu_j\right)^2} \left(\frac{t-\mu_j}{\sigma_j}\right)^2r_X^{'}\left(\frac{t-\mu_j}{\sigma_j}\right).$$
So, the second term of (\ref{2}) is positive (negative). Again, $u^3r_X^{''}(u)$ is increasing (decreasing) in $u$ gives
$$\left(\frac{t-\mu_i}{\sigma_i}\right)^3r_X^{''}\left(\frac{t-\mu_i}{\sigma_i}\right)\leq (\geq) \left(\frac{t-\mu_j}{\sigma_j}\right)^3r_X^{''}\left(\frac{t-\mu_j}{\sigma_j}\right).$$
Considering the fact that $r_X(u)$ is concave (convex) in $u$, from the above relation it can be written that,
$$\frac{1}{\left(t-\mu_i\right)^3}\left(\frac{t-\mu_i}{\sigma_i}\right)^3r_X^{''}\left(\frac{t-\mu_i}{\sigma_i}\right)\geq (\leq)\frac{1}{\left(t-\mu_j\right)^3} \left(\frac{t-\mu_j}{\sigma_j}\right)^3r_X^{''}\left(\frac{t-\mu_j}{\sigma_j}\right).$$
So, the first term of ($\ref{2}$) is also positive (negative). Thus $\frac{\partial\Psi}{\partial\mu_i}-\frac{\partial\Psi}{\partial\mu_j}\geq (\leq) 0$. So by Lemma \ref{l6} it can be concluded that $\Psi$ is s-convex (s-concave) in $\mbox{\boldmath $\mu$}$. Thus, using Lemma \ref{l10} the required result can be obtained.\\
When $\mbox{\boldmath $\mu$},\mbox{\boldmath $\sigma$}\in \mathcal{E}_+$, then the theorem can be proved in similar line. \hfill$\Box$\\
Combining Theorem \ref{th13} and Theorem \ref{th14} the following theorem can easily be obtained.
\begin{t1}\label{th15}
Let $\left\{X_1, X_2, \ldots, X_n\right\}$ and $\left\{Y_1, Y_2, \ldots, Y_n\right\}$ be two sets of independent random variables each following $LS$ model such that $X_i\sim $LS$\left(\lambda_i,\sigma_i,F\right)$ and $Y_i\sim $LS$\left(\mu_i,\sigma_i,F\right)$ for $i=1,2,\ldots,n$. Suppose further that $\mbox{\boldmath $\lambda$},\mbox{\boldmath $\mu$},\mbox{\boldmath $\sigma$}\in \mathcal{D}_+$ (or $\mathcal{E}_+$). Then, if $r_X(u)$ is increasing and concave (decreasing and convex) in $u$; $u^2r_X^{'}(u)$ is decreasing (increasing) in $u$ and $u^3r_X^{''}(u)$ is increasing (decreasing) in $u$ then, $\mbox{\boldmath $\lambda$}\preceq_w (\succeq_w) \mbox{\boldmath $\mu$}$ $\Rightarrow$ $X_{1:n}\leq_{lr}Y_{1:n}$.
\end{t1}
%%--------------------------comparison of parallel system (independent and dependet)----------------------------------------------------------
\section{Comparison of Parallel Systems with LS Distributed Components}
\setcounter{equation}{0}
\hspace*{0.3 in} Similar results for the parallel system are furnished in this section for dependent setup. These results are the extension of the Theorems 3.1-3.5 of Hazra \emph{et al.}~\cite{ha}, where the parallel system consists of independent $LS$ family distributed components. The following theorems \ref{th6}-\ref{th8} can be proved in the same line as of Theorem \ref{th1}, Theorem \ref{th2} and Theorem \ref{th3} respectively and hence the proofs are omitted. For results of independent random variables, one may refer to Hazra \emph{et al.}~\cite{ha}.\\
 %%--------------------------------Subsection 4.1-------------------------------------------------------
%\subsection{\small Heterogenous dependent Components}
\hspace*{0.3 in} Let $X$ and $Y$ be two random variables having distribution functions $F(\cdot)$ and $G(\cdot)$ respectively. Suppose that $X_i\sim LS\left(\lambda_i,\sigma_i,F\right)$ and $Y_i\sim LS\left(\mu_i,\xi_i,G\right)$ ($i=1,2,\ldots,n$) be two sets of $n$ dependent random variables with Archimedean copulas having generators $\psi_1$ (with $\phi_1=\psi_{1}^{-1}$) and $\psi_2$ ($\phi_2=\psi_{2}^{-1}$) respectively. Also suppose that $G_{n:n}\left(\cdot\right)$ and $H_{n:n}\left(\cdot\right)$ be the distribution functions of $X_{n:n}$ and $Y_{n:n}$ respectively. Then, 
\begin{equation*}
G_{n:n}\left(t\right)=\psi_1\left[\sum_{k=1}^n \phi_1\left\{F\left(\frac{t-\lambda_k}{\sigma_k}\right)\right\}\right],~t>max(\lambda_k,~\forall k),
\end{equation*}
and
\begin{equation*}
H_{n:n}\left(t\right)=\psi_2\left[\sum_{k=1}^n \phi_2\left\{G\left(\frac{t-\mu_k}{\xi_k}\right)\right\}\right],~t>max(\mu_k,~\forall k).
\end{equation*}
Let $\tilde r_X(u)$ and $\tilde r_Y(u)$ are the reversed hazard rate functions of the random variables $X$ and $Y$ respectively.
\begin{t1}\label{th6}
Let $X_1,X_2,...,X_n$ be a set of dependent random variables sharing Archimedean copula having generator $\psi_1$ such that $X_i\sim $LS$\left(\lambda,\sigma_i,F\right),~i=1,2,...,n$. Let $Y_1,Y_2,...,Y_n$ be another set of dependent random variables sharing Archimedean copula having generator $\psi_2$ such that $Y_i\sim $LS$\left(\lambda,\xi_i,G\right),~i=1,2,...,n$. Assume that $\mbox{\boldmath $\lambda$},\mbox{\boldmath $\sigma$},\mbox{\boldmath $\xi$}\in \mathcal{D}_+$ (or $\mathcal{E}_+$). Further suppose that $\phi_2\circ\psi_1$ is super-additive, $\psi_1$ or $\psi_2$ is log-convex and $X\ge_{st}Y$. If either $\tilde r_X(u)$ or $\tilde r_Y(u)$ is decreasing in $u$ then, $\mbox{\boldmath $\frac{1}{\sigma}$}\stackrel{w}{\succeq} \mbox{\boldmath $\frac{1}{\xi}$}$ $\Rightarrow$ $X_{n:n}\ge_{st}Y_{n:n}$.
\end{t1}

\begin{t1}\label{th7}
Let $X_1,X_2,...,X_n$ be a set of dependent random variables sharing Archimedean copula having generator $\psi_1$ such that $X_i\sim $LS$\left(\lambda_i,\sigma_i,F\right),~i=1,2,...,n$. Let $Y_1,Y_2,...,Y_n$ be another set of dependent random variables sharing Archimedean copula having generator $\psi_2$ such that $Y_i\sim $LS$\left(\lambda_i,\xi_i,G\right),~i=1,2,...,n$. Assume that $\mbox{\boldmath $\lambda$},\mbox{\boldmath $\sigma$},\mbox{\boldmath $\xi$}\in \mathcal{D}_+$ (or $\mathcal{E}_+$). Further suppose that $\phi_2\circ\psi_1$ is super-additive, $\psi_1$ or $\psi_2$ is log-convex and $X\ge_{st}Y$, then,
\begin{enumerate}
\item[i)] $\mbox{\boldmath $\frac{1}{\sigma}$}\stackrel{p}{\succeq} \mbox{\boldmath $\frac{1}{\xi}$}$ $\Rightarrow$ $X_{n:n}\ge_{st}Y_{n:n}$,$~$ if either $ur_X(u)$ or $ur_Y(u)$ is decreasing in $u$;
\item[ii)] $\mbox{\boldmath $\frac{1}{\sigma}$}\stackrel{rm}{\succeq} \mbox{\boldmath $\frac{1}{\xi}$}$ $\Rightarrow$ $X_{n:n}\ge_{st}Y_{n:n}$,$~$ if either $u^2r_X(u)$ or $u^2r_Y(u)$ is decreasing in $u$.
\end{enumerate}
\end{t1}
 
\begin{t1}\label{th8}
Let $X_1,X_2,...,X_n$ be a set of dependent random variables sharing Archimedean copula having generator $\psi_1$ such that $X_i\sim $LS$\left(\lambda_i,\sigma_i,F\right),~i=1,2,...,n$. Let $Y_1,Y_2,...,Y_n$ be another set of dependent random variables sharing Archimedean copula having generator $\psi_2$ such that $Y_i\sim $LS$\left(\mu_i,\sigma_i,G\right),~i=1,2,...,n$. Assume that $\mbox{\boldmath $\lambda$},\mbox{\boldmath $\sigma$},\mbox{\boldmath $\mu$}\in \mathcal{D}_+$ (or $\mathcal{E}_+$). Further suppose that $\phi_2\circ\psi_1$ is super-additive, $\psi_1$ or $\psi_2$ is log-convex and either $ur_X(u)$ or $ur_Y(u)$ is decreasing in $u,$ then, $X\ge_{st}Y$ and $\mbox{\boldmath $\lambda$}\succeq_w \mbox{\boldmath $\mu$}$ $\Rightarrow$ $X_{n:n}\ge_{st}Y_{n:n}$.
\end{t1} 
%\subsection{\small Heterogenous independent Components}

%%-----------------------------------References------------------------------------------------------------------

%%------------------------------About the authors----------------------------------------------------------
%\begin{center}
%\section*{ \small AUTHORS}
%\end{center}
 %%-------------------------------------------THE END-----------------------------------------------------------
\end{document}